\documentclass[a4paper,11pt,preprint]{imsart}
\usepackage{fullpage}

\usepackage{amsmath}
\usepackage[utf8]{inputenc}
\usepackage{type1cm}
\usepackage{amsfonts}
\usepackage{amssymb}
\usepackage{amscd}
\usepackage{amsbsy}
\usepackage{amsthm}
\usepackage{mathrsfs}
\usepackage{graphicx}
\usepackage{color}

\usepackage[pdfstartview=FitH]{hyperref}

\makeatletter
\@ifundefined{coloneqq}{%
\newcommand{\coloneqq}{\mathrel{\colon\mkern-10mu=}}}{}
\@ifundefined{eqqcolon}{%
\newcommand{\eqqcolon}{\mathrel{=\mkern-1mu\colon\mskip-6mu}}}{}
\def\blfootnote{\xdef\@thefnmark{}\@footnotetext}
\makeatother

\newcommand{\lt}{\left}
\newcommand{\rt}{\right}

\newcommand{\Knn}{K_{n,n}}
\newcommand{\Kn}{K_n}

\newcommand{\Knalpha}{K_{n,n/\alpha}}
\newcommand{\nbyalpha}{\lceil n/\alpha\rceil}
\newcommand{\Mn}{M_n}
\newcommand{\Mnalpha}{M_n^\alpha}
\newcommand{\calpha}{c_*^\alpha}

\newcommand{\To}{\rightarrow}
\newcommand{\Real}{\mathbf{R}}
\newcommand{\Natural}{\mathbf{N}}
\newcommand{\RPlus}{\Real_{+}}
\DeclareMathOperator*{\argmin}{arg\,min}
\newcommand{\parl}[1]{\left(#1\right)}
\newcommand{\brackl}[1]{\left[#1\right]}
\newcommand{\set}[1]{\left\{#1\right\}}

\newcommand{\ind}{\mathbf{1}}
\newcommand{\card}[1]{\left\lvert#1\right\rvert}
\newcommand{\ceil}[1]{\left\lceil#1\right\rceil}
\newcommand{\sEl}[1]{\sE{\left[#1\right]}}

\def\bigO{\mathop{O}\nolimits}
\def\bigOmega{\mathop{\Omega}\nolimits}

\def\bigTheta{\mathop{\Theta}\nolimits}

\newcommand{\limn}{\lim_{n\To\infty}}
\newcommand{\limk}{\lim_{k\To\infty}}
\newcommand{\limsupn}{\limsup_{n\To\infty}}

\newcommand{\liminfn}{\liminf_{n\To\infty}}

\newcommand{\cvglw}{\xrightarrow{\text{\upshape\tiny l.w.}}}   
\newcommand{\cvgltwo}{\xrightarrow{L^2}}   

\newcommand{\sumver}[2]{\sum_{#1\in#2(\phi)}}
\newcommand{\maxver}[2]{\max_{#1\notin#2(\phi),#1\sim\phi}}
\newcommand{\nbr}[1]{N_{#1}}
\newcommand{\graph}{{\mathcal{G}}}
\newcommand{\graphs}{{\mathcal{G}_*}}
\newcommand{\graphss}{{\mathcal{G}_{**}}}
\newcommand{\intG}[2]{\int_\graphs\sum_{v\sim\phi}#2\rd#1([G,\phi])}
\newcommand{\nbd}[1]{\mathcal{N}_\rho(#1)}
\newcommand{\Vinf}{\mathcal{V}}
\newcommand{\Einf}{\mathcal{E}}
\newcommand{\dirEdge}[1]{\overrightarrow{E}(#1)}

\newcommand{\Xone}{X_\phi^{(1)}}
\newcommand{\Xtwo}{X_\phi^{(2)}}
\newcommand{\pwit}{\mathcal{T}}

\newcommand{\len}[1]{\xi\left(#1\right)}
\newcommand{\leng}[2]{\xi_{#2}\left(#1\right)}

\newcommand{\lennalpha}[1]{\leng{#1}{\Knalpha}}
\newcommand{\lenT}[1]{\leng{#1}{\pwit}}
\newcommand{\X}[1]{X\mathopen{}\left(#1\right)}
\newcommand{\Xg}[3]{X_{#2}^{#3}\left(#1\right)}

\newcommand{\Xnalpha}[2]{\Xg{#1}{\Knalpha}{#2}}
\newcommand{\XT}[2]{\Xg{#1}{\pwit}{#2}}

\newcommand{\mapg}[3]{\pi_{#2}^{#3}#1}

\newcommand{\mapnalpha}[2]{\mapg{#1}{\Knalpha}{#2}}
\newcommand{\mapT}[2]{\mapg{#1}{\pwit}{#2}}
\newcommand{\mudiag}{\mu^\nearrow}

\newcommand{\nbrneq}[3]{{#1\sim#2,#1\neq#3}}
\newcommand{\mins}[3]{\min_{\nbrneq{#1}{#2}{#3}}}
\newcommand{\diff}[2]{\left(\len{#1,#2}-\X{#1,#2}\right)^+}
\newcommand{\diffg}[4]{\left(\leng{#1,#2}{#3}-\Xg{#1,#2}{#3}{#4}\right)^+}

\newcommand{\diffT}[3]{\diffg{#1}{#2}{\pwit}{#3}}

\def\polylog{\mathop{{\rm Li}_2}\nolimits}    

\newcommand{\alphaenumerate}[1]{\let\restoretheenumi\theenumi
\let\restorelabelenumi\labelenumi
\renewcommand{\theenumi}{\alph{enumi}}
\renewcommand{\labelenumi}{(\theenumi)}
#1
\let\theenumi\restoretheenumi
\let\labelenumi\restorelabelenumi
}
\renewcommand{\theenumi}{(\alph{enumi})}
\renewcommand{\labelenumi}{\theenumi}

\makeatletter
\@ifundefined{theorem}{
\@ifundefined{chapter}{ 
\newtheorem{theorem}{\noindent Theorem}
\newtheorem{lemma}{\noindent Lemma}
\newtheorem{corollary}{\noindent Corollary}
\newtheorem{proposition}{\noindent Proposition}
\newtheorem{conjecture}{\noindent Conjecture}
\newtheorem{assertion}{\noindent Assertion}
\newtheorem{assumption}{\noindent Assumption}
\newtheorem{condition}{\noindent Condition}
\newtheorem*{defnn}{\noindent Definition}
}{ 
\newtheorem{theorem}{\noindent Theorem}[chapter]
\newtheorem{lemma}{\noindent Lemma}[chapter]
\newtheorem{corollary}{\noindent Corollary}[chapter]
\newtheorem{proposition}{\noindent Proposition}[chapter]
\newtheorem{conjecture}{\noindent Conjecture}[chapter]
\newtheorem{assertion}{\noindent Assertion}[chapter]
\newtheorem{assumption}{\noindent Assumption}[chapter]
\newtheorem{condition}{\noindent Condition}[chapter]
\newtheorem*{defnn}{\noindent Definition}
}
}
{}
\makeatother

\newcommand{\rd}{\mathrm{d}} 
\DeclareMathOperator\sd{d\!} 
\DeclareMathOperator\sP{P}   
\DeclareMathOperator\sE{E}   
\newcommand{\cvgpr}{\xrightarrow{\text{\upshape\tiny P}}}   
\newcommand{\cvgwk}{\xrightarrow{\text{\upshape\tiny W}}}   
\newcommand{\cvgdist}{\xrightarrow{\text{\upshape\tiny D}}} 
\newcommand{\eqdist}{\stackrel{\text{\upshape\tiny D}}{=}}  
\chardef\ii="10

\DeclareMathAlphabet{\bi}{OT1}{ptm}{b}{it}

\newcommand{\bxi}{\boldsymbol{\xi}}

\newcommand{\allcomment}[1]{\relax}

\renewcommand{\rd}{\mathop{}\!\mathrm{d}}

\newcommand{\matching}{\mathcal{M}}
\newcommand{\Mopt}{\mathcal{M}_{\mathrm{opt}}}

\newcommand{\mintwo}[1]{{\min_{#1}}^{(2)}}

\renewcommand{\theenumi}{\alph{enumi}}
\renewcommand{\labelenumi}{(\theenumi)}

\begin{document}
\begin{frontmatter}

\title{Belief propagation for minimum weight many-to-one matchings in the random complete graph}
\runtitle{Belief propagation for many-to-one matchings}


\author{\fnms{Mustafa} \snm{Khandwawala}\thanksref{t1}\ead[label=e1]{mustafa@ece.iisc.ernet.in}}
\thankstext{t1}{The author is currently at INRIA, Paris, France.}
\runauthor{Khandwawala}

\affiliation{Indian Institute of Science}

\address{Department of Electrical Communication Engineering\\
Indian Institute of Science\\
Bangalore 560012, India\\
\printead{e1}}



\begin{abstract}
In a complete bipartite graph with vertex sets of cardinalities $n$ and $m$, assign random weights from exponential distribution with mean 1, independently to each edge. We show that, as $n\To\infty$, with $m=\nbyalpha$ for any fixed $\alpha>1$, the minimum weight of many-to-one matchings converges to a constant (depending on $\alpha$). Many-to-one matching arises as an optimization step in an algorithm for genome sequencing and as a measure of distance between finite sets. We prove that a belief propagation (BP) algorithm converges asymptotically to the optimal solution. We use the objective method of Aldous to prove our results. We build on previous works on minimum weight matching and minimum weight edge-cover problems to extend the objective method and to further the applicability of belief propagation to random combinatorial optimization problems.
\end{abstract}



\end{frontmatter}


\section{Introduction} \label{sec:many:introduction}
We address here an optimization problem over bipartite graphs related to the matching problem -- called the many-to-one matching problem, which has applications in an algorithm for genome sequencing and as a measure of distance between finite sets.
Given two sets $A$ and $B$, with $\card{A}\ge \card{B}$, consider a bipartite graph $G$ with vertex set $V=A\cup B,$ and edge set $E\subset A\times B$. A many-to-one matching in $G$ is a subgraph $\matching$ such that each vertex in $A$ -- called the {\em one side} -- has degree 1 in $\matching,$ and each vertex in $B$ -- called the {\em many side} -- has degree 1 or more in $\matching$. A many-to-one matching can be viewed as an onto function from $A$ to $B$. Each edge $e\in E$ has a weight $\xi_e \in \RPlus = [0,\infty)$. The cost of the matching $\matching$ is the sum of the weights of the edges in $\matching$.
We focus here on the minimum cost many-to-one matchings on complete bipartite graphs where the edge-weights are independent and identically distributed (i.i.d.) random variables. 

For several optimization problems over graphs with random edge-weights, where the objective is to minimize the sum weight of a subset of edges under some constraints, the expected optimal cost converges to a constant as $n\To\infty$. Such results are known for the minimum spanning tree \cite{Fri1985}, matching \cite{Ald1992,Ald2001}, edge-cover \cite{HesWas2010,KhaSun2012}, and traveling salesman problem (TSP) \cite{Was2010} over complete graphs or complete bipartite graphs with i.i.d. edge-weights. Furthermore, simple, iterative \emph{belief propagation} algorithms can be used to find asymptotically optimal solutions for matching \cite{SalSha2009} and edge-cover \cite{KhaSun2012}.

We prove a similar result for many-to-one matching. Let us be a little more specific.
Fix a real number $\alpha>1$. In the complete bipartite graph $K_{n,m}$ with vertex sets $A$ and $B$ of cardinality $n$ and $m=\nbyalpha$ respectively, we take the edge-weights to be i.i.d.\ random variables having the exponential distribution with mean $1$. Denote by $\Mnalpha$ the minimum cost of many-to-one matching on $K_{n,m}$. 
We show that the expected value of $\Mnalpha$ converges to a constant (which depends on $\alpha$ and can be computed) as $n\To\infty$. Further, we show that a belief propagation (BP) algorithm finds the asymptotically optimal many-to-one matching in $O(n^2)$ steps.

We proceed to show these results via the objective method following \cite{Ald2001} for matching and \cite{KhaSun2012} for edge-cover.
Before we give the background and overview of the methods, we describe two practical applications of the many-to-one matching problem.

\section{Applications of many-to-one matching}
\subsection{Restriction scaffold problem}
Ben-Dor et al.\ \cite{BenKarSch2003} introduced the \emph{restriction scaffold problem} for assisting the finishing step in a shotgun genome sequencing project. Shotgun sequencing refers to the process of generating several short clones of a large target genome and then sequencing these individual clones. Sequencing a short genome is much easier than sequencing the entire genome, and the cloning process is automated and fast, but the locations of clones are random. This process, which is like a random sampling, leaves us with several sequenced segments of the target. It is possible to construct the target sequence by assembling these short sequences by matching their overlaps, provided that we have sufficient redundancy in the samples to cover the entire sequence with enough overlap.

\begin{figure}[t]
 \centering
 \includegraphics{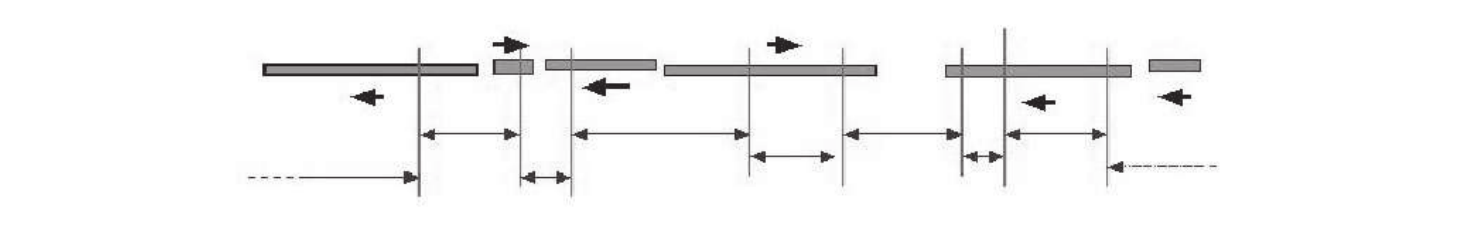} 
 
 \caption{``Fraction of a contiguration and one enzyme data. Top: the contigs; bold arrow denote the orientation of each contig.
Vertical lines: restrictions sites. Bottom: restriction fragments between neighboring sites.''\newline
Image and the quoted caption from \cite{BenKarSch2003}.\label{fig:rsp}}
\end{figure}

Usually, asking for enough redundancy -- that this assembly based on overlaps reconstructs the whole sequence -- will result in excessive cloning and sequencing tasks. As a tradeoff, we can give up on this high redundancy, and accept as the output of this assembly task large contiguous sequences (Ben-Dor et al.\ call these \emph{contigs}), which together cover a large part of the target sequence. If the gaps -- the positions in the target sequence that are not covered by these contigs -- are small, those areas of the genome can be targeted manually to finish the sequencing project. However, the relative position and orientation of the contigs, and the sizes of gaps between successive contigs are unknown. We need some information to figure out how these contigs align along the target genome. See the top part of figure \ref{fig:rsp}.

Ben-Dor et al.\ propose the use of restriction enzymes to generate this information. Restriction enzymes cut the DNA at specific recognition sequences (differing according to the enzyme) called the restriction site of the enzyme. These restriction sites form a scaffold over which the contigs can be assembled. Again, see figure \ref{fig:rsp}. They measure the approximate sizes of the fragments obtained by digesting the target DNA with a restriction enzyme. On the other hand, by assuming a particular ordering and orientation of the contigs and the gaps between successive contigs the sizes of the fragments can be computed by reading out the recognition sequence of the enzyme along the contigs (their sequences are known). Comparing these computed fragment sizes with the measured fragment sizes for several different restriction enzymes, gives them a way to obtain a new estimate of the ordering, orientation and gaps. This is the basis of their iterative algorithm for solving the restriction scaffold problem.

The step of the iterative algorithm involving the comparison of the computed fragment sizes with the measured fragment sizes is formulated as a many-to-one matching problem -- the computed fragment sizes on the ``one'' side, and the measured fragment sizes on the ``many'' side. The measurement step process does not indicate the number of fragments corresponding to a particular measurement. The many-to-one nature of the association captures the notion that a measurement may correspond to one or more physical fragments, and each physical fragment must correspond to one measurement. The weight of an edge captures how much the two sizes differ. Specifically, the weight of an edge joining a computed fragment size $s(c)$ to a measured fragment size $s(m)$ is $\card{\log(s(c))-\log(s(m))}$. They solve the many-to-one matching problem using this linear property: the vertices are represented as points on a line, a vertex corresponding to a size $s$ is at position $\log(s)$ on the line, and the cost of matching two 
points is the distance between the points.

We do not assume this linear structure in our formulation of the many-to-one problem on the complete bipartite graph. However the i.i.d. weight assumption is the mean-field approximation of a geometric structure in the following sense. Take a set $A$ of $n$ points and a set $B$ of $m=\nbyalpha$ points independently and uniformly at random in the unit sphere in $\Real^d$. When $m$ and $n$ are large, the probability that the distance between a typical point of $A$ and a typical point of $B$ is less than $r$ is approximately $r^d$ for small $r$. We can represent the distances between points of the two sets by forming a complete bipartite graph and taking the weight of an edge as the distance between the corresponding points. For a mean-field model we ignore the geometric structure of $\Real^d$ (and triangular inequality) and only take into account the interpoint distances, by taking the edge-weights to be i.i.d.\ copies of a nonnegative random variable $\xi$ satisfying $\Pr(\xi<r)\approx r^d$ as $r\To 0$.

\subsection{A measure of distance between sets of points}
A concept of distance between finite sets of points is useful in many areas like machine learning, computational geometry, and comparison of theories. Such a distance is derived from a given distance function (or a metric) between the points. For example, in a clustering process over a set of examples, suppose we are given a function $d$ such that $d(e_1,e_2)$ corresponds to the distance between two examples $e_1$ and $e_2$. Clustering methods such as TIC \cite{BloDerRam1998} rely on a function $d_1$ that specifies the distance $d_1(C_1,C_2)$ for two clusters (which are sets of examples) $C_1$ and $C_2$.

Eiter and Mannila \cite{EitMan1997} use the term \emph{surjection distance} for a measure of similarity, based on optimal many-to-one matching, between two finite sets of points in a metric space. The surjection distance is the minimum cost of a many-to-one matching, with the points of the larger set forming the vertices of the one side, the points of the smaller set forming the vertices of the many side of a bipartite graph, and the edge-weights taken as the distance between the corresponding points. By reduction to the minimum weight matching problem they find an $O(n^3)$ algorithm for many-to-one matching. Belief propagation (in the random setting we consider here) yields an asymptotically optimal solution in $O(n^2)$ steps.

\section{Background and overview of the methods}
Several of the results for random combinatorial optimization problems have originated as conjectures supported by heuristic statistical mechanical methods such as the replica method and the cavity method \cite{Par1998,MezPar1985,MezPar1986,MezPar1987,MezParVir1987}. Establishing the validity of these methods is a challenge even for mathematically simple models.

Aldous \cite{Ald1992,Ald2001} provided a rigorous proof of the $\zeta(2)$ limit conjecture for the random matching (or assignment) problem. The method he used is called the \emph{objective method}. A survey of this method and its application to a few other problems was given in \cite{AldSte2004}. The objective method provides a rigorous counterpart to the \emph{cavity method}, but has had success only in a few settings. 

The crux of the objective method is to obtain a suitable distributional limit for a sequence of finite random weighted graphs of growing vertex set size. The space of weighted graphs is endowed with a metric that captures the proximity based on the local neighborhood around a vertex. The weak convergence of probability measures on this space is called \emph{local weak convergence} to emphasize that the convergence applies to events that are essentially characterized by the local structure around a typical vertex, and not to events that crucially depend on the global structure of the graph. This form of convergence has appeared in varying forms in \cite{Haj1990,BenSch2001,AldSte2004}.
Once we have this convergence, we can relate the optimization problem on finite graphs to an optimization problem on the limit object. 
For example, Aldous \cite{Ald2001} showed that the limit object for the sequence of random weighted complete bipartite graphs $\Knn$ is the \emph{Poisson weighted infinite tree (PWIT)}. Aldous then writes the optimization problem in terms of local conditions at the vertices of the tree, and using the symmetry in the structure of the PWIT, these relations result in a distributional identity called the \emph{recursive distributional equation}. A solution to this equation is then used to describe asymptotically optimal matchings on $\Knn$.

One can exploit the recursive distributional equation to construct an iterative decentralized message passing algorithm, a versions of which is called belief propagation (BP) in the computer science literature. BP algorithms are known to converge to the correct solution on graphs without cycles. For graphs with cycles, provable guarantees are known for BP only for certain problems. For example, Bayati, Shah and Sharma \cite{BayShaSha2008} proved that the BP algorithm for maximum weight matching on bipartite graphs converges to the correct value as long as the maximum weight matching is unique, and Sanghavi, Malioutov and Willsky \cite{SanMalWil2011} proved that BP for matching is as powerful as LP relaxation. Salez and Shah \cite{SalSha2009} studied belief propagation for the random assignment problem, and showed that a BP algorithm on $\Knn$ converges to an update rule on the limit object PWIT. The iterates on the PWIT converge in distribution to the minimum cost assignment. The iterates are near the 
optimal solution in $O(n^2)$ steps whereas the worst case 
optimal algorithm on bipartite graphs is $O(n^3)$ (expected time $O(n^2 \log n)$ for i.i.d.\ edge capacities); see Salez and Shah \cite{SalSha2009} and references therein. 

In message passing algorithms like belief propagation, the calculations are performed at each vertex of a graph from the messages received from its neighbors. This local nature makes feasible its analysis via the objective method. The author and Sundaresan, in \cite{KhaSun2012}, extended the objective method by combining local weak convergence and belief propagation to prove and characterize the limiting expected minimum cost for the edge-cover problem. In this paper, we implement this general program for the many-to-one matching problem. The proof relies on two properties: (1) the neighborhoods around most vertices are tree-like, and (2) the effect of the boundary on the belief propagation messages at the root of the tree vanishes as the boundary is taken to infinity. The second property is formalized as ``endogeny" (see \cite{AldBan2005}) and is handled in Section~\ref{sec:many:endogeny}.


\section{Contributions of this paper}
This paper extends the line of work on combinatorial optimization problems over graphs with i.i.d.\ edge-weights. Although the methods are similar to those in \cite{Ald2001, SalSha2009, KhaSun2012}, there are key differences from earlier work in addressing the many-to-one matching problem. We now highlight them.
For example, we must deal with two types of vertices corresponding to two different constraints on the degrees. The different viewpoints from a vertex in $A$ and a vertex in $B$ makes the limit object for the sequence of graphs $K_{n,m}$ different from the Poisson weighted infinite tree (PWIT) obtained as the limit of bipartite graphs $\Knn$ with vertex sets of the same cardinality. The difference is in the distribution of the weights of the edges at odd and even levels from the root. As a consequence, the recursive distributional equation (RDE) -- guiding the construction of an optimal many-to-one matching on the limit object -- is expressed in terms of a two-step recursion. 

We prove that this recursive distributional equation has a unique solution, and that the many-to-one matching on the limit object generated by the solution of this equation is optimal among all \emph{involution invariant} many-to-one matchings. To establish correctness of BP, we establish that the fixed-point solution of the RDE has a full domain of attraction, and that the stationary process with the fixed-point marginal distribution on the limit object satisfies endogeny.


\section{Main results} \label{sec:many:results}
From now on we will write $\Knalpha$ to mean $K_{n,m},$ where $m=\nbyalpha$. Recall that $A$ and $B$ denote the \emph{one} side and the \emph{many} side in the bipartite graph $\Knalpha$. The weights of the edges are i.i.d. random variables with the exponential distribution, and for simplicity in the subsequent analysis, we take the mean of the edge-weights to be $n$. This only changes the cost of the optimal solution by a factor of $n$, and so we have a rescaling in the left hand side of (\ref{eq:many:limit}).

Our first result establishes the limit of the scaled expected minimum cost of the random many-to-one matching problem.
\begin{theorem} \label{thm:many:limit}
For $\alpha>1$, the minimum cost $\Mnalpha$ of many-to-one matching on $\Knalpha,$ satisfies
 \begin{align}
 \limn n^{-1}\sE{\Mnalpha} &=-\polylog(-\gamma^{-1})-\frac{1}{2}\log^2(1+\gamma^{-1}) +w_o\log(1+\gamma^{-1})+w_o\notag\\
 &\eqqcolon \calpha ,\label{eq:many:limit}
 \end{align}
where $\polylog$ is the dilogarithm: $\polylog(z)=-\int_0^z\frac{\log(1-t)}{t}\rd t=\sum_{k=1}^\infty\frac{z^k}{k^2}$; $w_o$ is the positive solution to $\alpha=w_o+e^{-w_o}$; $\gamma$ equals $w_o e^{w_o}$.
\end{theorem}

Our second result shows that a belief propagation algorithm gives a many-to-one matching that is asymptotically optimal as $n\To\infty$. The BP algorithm is based on the corresponding BP algorithms developed in \cite{SalSha2009,BayShaSha2008} for matching and in \cite{KhaSun2012} for edge-cover. The proof uses the technique of showing that the update rule of BP converges to an update rule on a limit infinite tree. 

We now define a BP algorithm on an arbitrary bipartite graph $G=(A\cup B,E)$ with edge-costs. First, some notation. We write $w\sim v$ to mean $w$ is a neighbor of $v,$ i.e., $\set{v,w}\in E$. For an edge $e=\set{v,w}\in E,$ we write its cost as $\leng{e}{G}$ or $\leng{v,w}{G}$ (with the ordering irrelevant). For each vertex $v\in V,$ we associate a nonempty subset $\mapg{(v)}{G}{}$ of its neighbors such that $\mapg{(v)}{G}{}$ has cardinality 1 if $v\in A$. Write $\matching(\mapg{}{G}{})=\set{\set{v,w}\mid v\in A\cup B, w\in\mapg{(v)}{G}{}}$.

The BP algorithm is an iterative message passing algorithm. In each iteration $k\ge 0,$ every vertex $v\in V$ sends a message $\Xg{w,v}{G}{k}$ to each neighbor $w\sim v$ according to the following rules:
\begin{description}
 \item[Initialization:] \begin{equation}\label{eq:many:bpInit}
                        \Xg{w,v}{G}{0}=0.
                       \end{equation}
 \item[Update rule:] \begin{equation}\label{eq:many:bpUpdate}
                        \Xg{w,v}{G}{k+1}=\begin{cases}
                        \mins{u}{v}{w}\set{\leng{v,u}{G}-\Xg{v,u}{G}{k}},\quad&\text{if }v\in A\\
                        \mins{u}{v}{w}\set{\diffg{v}{u}{G}{k}},\quad&\text{if }v\in B.
                                          \end{cases}
                       \end{equation}
 \item[Decision rule:] \begin{gather}\label{eq:many:bpDec}
                        \mapg{(v)}{G}{k}=\begin{cases}
                        \argmin_{u\sim v}\set{\leng{v,u}{G}-\Xg{v,u}{G}{k}},\quad&\text{if }v\in A\\
                        \argmin_{u\sim v}\set{\diffg{v}{u}{G}{k}},\quad&\text{if }v\in B.
                                         \end{cases}
                       \end{gather}
\end{description}

Note that the subset $\matching(\mapg{(v)}{G}{k})$ is not necessarily a many-to-one matching, but it can be modified to be a many-to-one matching $\matching_n^k$ by some patching. We also remark that while $\leng{v,u}{G}=\leng{u,v}{G},$ the messages $\Xg{w,v}{G}{k}$ and $\Xg{v,w}{G}{k}$ can differ.

We analyze the belief propagation algorithm for $G=\Knalpha$ and i.i.d.\ exponential random edge-costs, and prove that after sufficiently large number of iterates, the expected cost of the many-to-one matching given by the BP algorithm is close to the limit value in Theorem~\ref{thm:many:limit}.
\begin{theorem} \label{thm:many:bpconv}
On $\Knalpha$, the output $\matching\lt(\mapg{}{\Knalpha}{k}\rt)$ of the BP algorithm can be patched to get a many-to-one matching $\matching_n^k$ that satisfies
 \begin{equation}\label{eq:many:thmbpconv}  \limk\limn n^{-1}\sEl{\sum_{e\in\matching_n^k}\leng{e}{\Knalpha}}=\calpha.
 \end{equation}
\end{theorem}

\section{Local weak convergence} \label{sec:many:lwConv}
First we recollect the terminology for defining convergence of graphs, borrowed from \cite{AldSte2004}.

\subsection{Rooted geometric networks}
A graph $G=(V,E)$ along with a length function $l:E\To(0,\infty]$ and $\mathsf{label}:V\To\mathcal{L}$ is called a {\em network}. $\mathcal{L}$ is a finite set of labels. The {\em distance} between two vertices in the network is the infimum of the sum of lengths of the edges of a path connecting the two vertices, the infimum being taken over all such paths. We call the network a {\em geometric network} if for each vertex $v\in V$ and positive real $\rho,$ the number of vertices within a distance $\rho$ of $v$ is finite. We denote the space of geometric networks by $\graph$.

A geometric network with a distinguished vertex $v$ is called a {\em rooted geometric network} with root $v$. We denote the space of all connected rooted geometric networks by $\graphs$. It is important to note that in $\graphs$ we do not distinguish between rooted isomorphisms of the same network -- isomorphisms that preserve the graph structure along with the root, vertex labels, and edge-lengths. We will use the notation $(G,\phi)$ to denote an element of $\graphs$ which is the isomorphism class of rooted networks with underlying network $G$ and root $\phi$.

\subsection{Local weak convergence}
\begin{figure}[t]
\def\svgwidth{130px}
\centering
\input{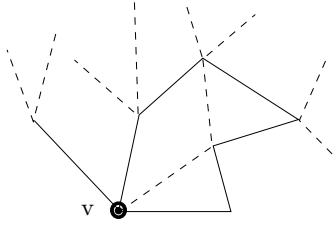}
\caption{\label{fig:many:nbd}Neighborhood $\nbd{G}$ of graph $G$. The solid edges form the neighborhood, and form paths of length at most $\rho$ from the root $v$. Dashed edges are the other edges of $G$.}
\end{figure}

We call a positive real number $\rho$ a {\em continuity point} of $G$ if no vertex of $G$ is exactly at a distance $\rho$ from the root of $G$. Let $\nbd{G}$ denote the neighborhood of the root of $G$ up to distance $\rho$. $\nbd{G}$ contains all vertices of $G$ which are within a distance $\rho$ from the root of $G$ (Figure~\ref{fig:many:nbd}). We take $\nbd{G}$ to be an element of $\graphs$ by inheriting the same length function $l$ as $G,$ same label function as $G,$ and the same root as that of $G$.

We say that a sequence of rooted geometric networks $G_n, n\ge 1$ {\em converges locally} to an element $G_\infty$ in $\graphs$ if for each continuity point $\rho$ of $G_\infty,$ there is an $n_\rho$ such that for all $n\ge n_\rho,$ there exists a graph isomorphism $\gamma_{n,\rho}$ from $\nbd{G_\infty}$ to $\nbd{G_n}$ that maps the root of the former to the root of the latter, preserves the labels of the vertices, and for each edge $e$ of $\nbd{G_\infty},$ the length of $\gamma_{n,\rho}(e)$ converges to the length of $e$ as $n\To\infty$.

The space $\graphs$ can be suitably metrized to make it a separable and complete metric space. See \cite{AldLyo2007} for details on the metric. One can then consider probability measures on this space and endow that space with the topology of weak convergence of measures. This notion of convergence is called {\em local weak convergence}.

In our setting of complete graphs $\Knalpha=(A\cup B,E)$ with random i.i.d.\ edge-costs $\set{\xi_e,e\in E},$ we regard the edge-costs to be the lengths of the edges, and declare a vertex $\phi$ of $\Knalpha$ chosen uniformly at random as the root of $\Knalpha$. Assign the label $o$ to the vertices in $A,$ and the label $m$ to the vertices in $B$; the label of a vertex will indicate its permitted degrees in a many-to-one matching: $o$ for one and $m$ for many. This makes $(\Knalpha,\phi),$ the isomorphism class of $\Knalpha$ and the root $\phi,$ a random element of $\graphs$. We will show (Theorem~\ref{thm:many:pwitConv} below) that the sequence of random rooted geometric networks $\Knalpha$ converges in the local weak sense to an element of $\graphs,$ which we call the {\em Poisson weighted infinite tree (PWIT)}.

\subsection{Poisson weighted infinite tree}
The graph that we define here is a variant of the Poisson weighted infinite tree that appears in \cite{AldSte2004}. We use the notation from \cite{SalSha2009} to define the variant.
This object is a random element of $\graphs$ whose distribution is a mixture of the distributions of two random trees: $\pwit_\alpha^o$ and $\pwit_\alpha^m$. Let $\Vinf$ denote the set of all finite strings over the alphabet $\Natural=\set{1,2,3,\ldots}$. Let $\phi$ denote the empty string and ``." the concatenation operator. Let $\card{v}$ denote the length of the string $v$. For a string $v$ with $\card{v}\ge 1,$ let $\dot{v}$ denote the parent of $v,$ i.e., the string with the last letter deleted. Let 
\[
\Einf=\set{\set{v,v.i},v\in\Vinf,i\in\Natural}.
\]
Both $\pwit_\alpha^o$ and $\pwit_\alpha^m$ have the same underlying graph $(\Vinf,\Einf)$ with root $\phi$. Each vertex $v$ in $\Vinf$ has a label from the set $\mathcal{L}=\set{o,m}$: $v$ has label $o$ in $\pwit_\alpha^o$ (respectively $\pwit_\alpha^m$) if $\card{v}$ is even (respectively odd), and label $m$ if $\card{v}$ is odd (respectively even). In particular, the root $\phi$ has label $o$ in $\pwit_\alpha^o,$ and label $m$ in $\pwit_\alpha^m$. Construct a family of independent Poisson processes, with intensity 1, on $\RPlus$: 
\[
\set{\bxi^v=(\xi^v_1,\xi^v_2,\ldots),v\in\Vinf}.
\]
Take the length of an edge $\set{v,v.i}$ in $\Einf$ to be $\xi^v_i$ if $v$ has label $m,$ and $\alpha\xi^v_i$ if $v$ has label $o$. If we denote the probability measure of $\pwit_\alpha^o$ and $\pwit_\alpha^m$ by $\nu_{o}$ and $\nu_{m}$ respectively, then the probability measure of $\pwit_\alpha$ is the mixture
\begin{equation*}
\nu=\frac{\alpha}{1+\alpha}\nu_o+\frac{1}{1+\alpha}\nu_m.
\end{equation*}
We now prove that this distribution is the local weak limit of the distributions of $(\Knalpha,\phi),$ where $\phi$ is chosen uniformly at random from the vertices of $\Knalpha$.

\begin{theorem} \label{thm:many:pwitConv}
The sequence $(\Knalpha,\phi)$ having edge-weights that are independent random variables with exponential mean $n$ distribution converges to $\pwit_\alpha$ (with root $\phi$) as $n\To\infty$ in the local weak sense.
\end{theorem}
\begin{proof}
In $\Knalpha,$ let $A$ denote the larger set with cardinality $n$ and let $B$ denote the set with smaller cardinality $m=\ceil{n/\alpha}$. Assign the label $o$ to vertices in $A$ (the one side), and label $m$ to the vertices in $B$ (the many side). Condition the root $\phi$ to have label $o$: this occurs with probability $n/(n+\nbyalpha)$. Fix a positive integer $N,$ and for $n>N$ identify vertices of $\Knalpha$ with strings in $\Vinf$ restricted to the alphabet $\set{1,2,\ldots,N}$. The strings $1,2,\ldots,N$ denote the first $N$ neighbors of the root $\phi$ arranged in increasing order of the weights of edges incident on the root. If $v\in\Vinf\setminus\phi$ denotes some vertex in $\Knalpha$ then $v.1,v.2,\ldots,v.N$ denote the first $N$ neighbors of $v$ excluding $\dot{v}$ arranged in increasing order of the weights of edges  incident on $v$. This generates a subtree of width $N$ of $(\Vinf,\Einf)$. Assign weights and vertex labels to the edges and vertices of the subtree according to those in $\
Knalpha$. Call this weighted tree $T_{N,n}$. Note that several vertices in $\Vinf$ may correspond to a single vertex of $\Knalpha$. Fix another positive integer $H$. Call the restriction of $T_{N,n}$ to vertices $v$ with $\card{v}\le H$ the weighted tree $T_{N,H,n}$.

Let (the scripted) $\pwit_{N}$ and $\pwit_{N,H}$ denote the corresponding restrictions of $\pwit_\alpha^o$. Clearly, $T_{N,H,n}$ and $\pwit_{N,H}$ are isomorphic finite graphs with the same vertex labels (The roots of both are labeled $o$). We can view the edge-weights as elements of $\Real^{\Einf(\pwit_{N,H})},$ where $\Einf(\pwit_{N,H})$ denotes the edges in $\pwit_{N,H}$. Write $\nu_{N,H,n}$ and $\nu_{N,H}$ for the measure of $T_{N,H,n}$ and $\pwit_{N,H}$ respectively.

We bound the likelihood ratio $\frac{\sd \nu_{N,H,n}}{\sd \nu_{N,H}}$ for the restriction that each vertex of $T_{N,H,n}$ corresponds to a unique vertex in $\Knalpha$. For any vertex $v,$ the vertex $v.1$ is taken to be the $w$ that minimizes $\xi(v,w)$ among all $w\neq\dot{v},$ and $w\sim v$. Suppose that we have associated vertices upto $v.i$ for some $1\le i<N$ to vertices of $\Knalpha$. If $v$ has label $o$ the next vertex $v.(i+1)$ is the vertex $w$ for which the difference $\xi(v,w)-\xi(v,v.i),w\notin \set{\dot{v},v.1,\ldots,v.i}$ is the minimum. The density of the difference $\xi(v,v.(i+1))-\xi(v,v.i)$ for $i\ge 1$ (respectively, $\xi(v,v.1)$ for $i=0$) at $x>0$ is at least 
\begin{equation} \label{eq:many:diffdensity}
\left(1-\frac{\card{\pwit_{N,H}}}{\nbyalpha-N}\right) \brackl{\frac{\nbyalpha-1-i}{n}\exp\left(-\frac{\nbyalpha-1-i}{n}x\right)},
\end{equation}
where $\card{\pwit_{N,H}}$ is the number of vertices in $\pwit_{N,H}$.
To see this, observe that the term with square brackets is, using the memoryless property of the exponential distribution, the density of the minimum among the $\nbyalpha-(i+1)$ exponential random variables corresponding to nodes on the many side (excluding $\set{\dot{v},v.1,\ldots,v.i}$). The term outside the square brackets is a lower bound on the probability that the minimum is not at one of the vertices of $\Knalpha$ already picked. This establishes (\ref{eq:many:diffdensity}).
The density of the corresponding difference in $\pwit_{N,H}$ is 
\[
\frac{1}{\alpha}\exp(-\frac{x}{\alpha}).
\]
The ratio of the conditional densities is at least
\[
\left(1-\frac{\card{\pwit_{N,H}}}{n/\alpha-N}\right)\left(1-\frac{N}{n/\alpha}\right).
\]
Doing this for all vertices of $\pwit_{N,H},$ and using the fact that under the restriction all edges are distinct, we get
\begin{equation*}
\frac{\rd \nu_{N,H,n}}{\rd \nu_{N,H}}\ge\left(\left(1-\frac{\card{\pwit_{N,H}}}{n/\alpha-N}\right)\left(1-\frac{N}{n/\alpha}\right)\right)^{\card{\pwit_{N,H}}}.
\end{equation*}
The lower bound converges to $1$ as $n\To\infty$. This implies that 
\begin{equation}\label{eq:many:cvgNH}
\nu_{N,H,n}\cvgwk\nu_{N,H},
\end{equation}
as $n\to\infty$.

Now to get the local weak convergence result, fix a $\rho>0$ and $\epsilon>0,$ and take $N$ and $H$ large enough such that uniformly for all $n>\card{\pwit_{N,H}}$ the probability that the $\rho$ neighborhood of the root of $\Knalpha$ is a subset of $T_{N,H,n}$ is greater than $1-\epsilon$; under this event the $\rho$ neighborhood of the root of $\Knalpha$ converges weakly to the $\rho$ neighborhood of the root of $\pwit_\alpha^o$ by (\ref{eq:many:cvgNH}). This shows that conditional on the root having label $o,$ $\pwit_\alpha^o$ is the local weak limit of $\Knalpha$. Similarly we can show that conditional on the root having label $m,$ the local weak limit is $\pwit_\alpha^m$.
The probability of the label of the root to be $o$ and $m$ converges respectively to $\alpha/(1+\alpha)$ and $1/(1+\alpha)$. This completes the proof.
\end{proof}

The above theorem is a generalization of Aldous's result \cite[Lemma~10]{Ald1992} for the local weak limit of complete bipartite graphs in which both the vertex subsets have the same cardinality. A similar result was earlier established by Hajek \cite[Sec.~IV]{Haj1990} for a class of sparse Erd\H{o}s-R\'{e}nyi random graphs.

The above theorem says that if we look at an arbitrary large but fixed neighborhood of the root of $\Knalpha,$ then for large $n,$ it looks like the corresponding neighborhood of the root of $\pwit$. This suggests that if boundary conditions can be ignored we may be able to relate optimal many-to-one matchings on $\Knalpha$ with an appropriate many-to-one matching on $\pwit$ (to be precise, an optimal {\em involution invariant} many-to-one matching 
 on the PWIT). Furthermore, the local neighborhood of the root of $\Knalpha$ is a tree for large enough $n$ (with high probability). So we may expect belief propagation on $\Knalpha$ to converge. Both these observations form the basis of the proof technique here, just like in the case of edge-cover in \cite{KhaSun2012}. In the matching case these ideas were developed by Aldous in \cite{Ald1992,Ald2001} to obtain the limit, and later by Salez and Shah \cite{SalSha2009} to show the applicability of belief propagation.

 First, let us see that the expected cost of an optimal many-to-one matching on $\Knalpha$ can be written in terms of the weights of the edges from the root:
 \[
\sEl{ \sum_{e\in M_n^*} \lennalpha{e}} = n \sEl{ \sum_{\set{\phi,v}\in M_n^*} \len{\phi,v} \mid \mathsf{label}(\phi)=o}.
 \]
This is because the root, conditioned to have label $o,$ is uniform among the $n$ vertices of the subset $A$. We seek to relate the expectation on the right with the expectation corresponding to optimal many-to-one matchings on the PWIT.

\section{Recursive distributional equation} \label{sec:many:rde}
For a labeled tree $\pwit$ with root $\phi,$ a many-to-one matching is a collection of edges such that vertices with label $o$ have degree 1 and the vertices with label $m$ have degree at least 1. The cost is the sum of the weights of the edges in the collection. This is well defined for a finite $\pwit$. Write $C(\pwit)$ for the minimum cost of many-to-one matching on $\pwit,$ and write $C(\pwit\setminus\phi)$ for the minimum cost of many-to-one matching on $\pwit\setminus\phi,$ which is the forest obtained after removing the root $\phi$ from $\pwit$. Assume that the size of the tree is such that such many-to-one matchings exist. We derive a recursive relation involving these quantities.

If $v$ is a child of the root $\phi$ of $\pwit$ let $\pwit_v$ denote the subtree containing $v$ and all its descendants, with $v$ as the root. 
Irrespective of the label of $\phi,$ we have
\begin{equation}\label{eq:many:costwophi}
 C(\pwit\setminus\phi)=\sum_{v\sim\phi} C(\pwit_v).
\end{equation}
If $\phi$ has label $o$ then among all matchings which match $\phi$ to a vertex $u$ the minimum cost is
\[
\len{\phi,u}+\min\set{C(\pwit_u),C(\pwit_u\setminus u)}+\sum_{v\sim\phi,v\neq u}C(\pwit_v).
\]
$C(\pwit)$ is the minimum of the above over all children of the root:
\begin{equation} \label{eq:many:costwithphio}
 C(\pwit)=\min_{u\sim\phi}\set{\len{\phi,u}+\min\set{C(\pwit_u),C(\pwit_u\setminus u)}+\sum_{v\sim\phi,v\neq u}C(\pwit_v)}.
\end{equation}
Subtracting (\ref{eq:many:costwophi}) from (\ref{eq:many:costwithphio}), we get
\begin{equation}\label{eq:many:costdiffo}
 C(\pwit)-C(\pwit\setminus\phi)=\min_{u\sim\phi}\set{\len{\phi,u}-\parl{C(\pwit_u)-C(\pwit_u\setminus u)}^+}.
\end{equation}
If $\phi$ has label $m$ then among all matchings which match $\phi$ to vertices in a nonempty set $I$ the minimum cost is
\[
\sum_{u\in I}(\len{\phi,u}+C(\pwit_u\setminus u))+\sum_{v\sim\phi,v\notin I}C(\pwit_v).
\]
$C(\pwit)$ is the minimum of the above over all nonempty $I$:
\begin{equation}  \label{eq:many:costwithphim}
 C(\pwit)=\min_{I\text{ nonempty}}\set{\sum_{u\in I}(\len{\phi,u}+C(\pwit_u\setminus u))+\sum_{v\sim\phi,v\notin I}C(\pwit_v).}
\end{equation}
Subtracting (\ref{eq:many:costwophi}) from (\ref{eq:many:costwithphim}), we get
\begin{equation*}
 C(\pwit)-C(\pwit\setminus\phi)=\min_{I\text{ nonempty}}\set{\sum_{u\in I}\brackl{\len{\phi,u}-(C(\pwit_u)-C(\pwit_u\setminus u))}},
\end{equation*}
the positive part of which is
\begin{equation}\label{eq:many:costdiffm}
\parl{C(\pwit)-C(\pwit\setminus\phi)}^+=\min_{u\sim\phi}\set{\brackl{\len{\phi,u}-(C(\pwit_u)-C(\pwit_u\setminus u))}^+}.
\end{equation}

Now consider what happens when we take $\pwit$ to be the infinite tree $\pwit_\alpha$ (PWIT). The quantities $C(\pwit_\alpha)$ and $C(\pwit_\alpha\setminus\phi)$ become infinite, and the relations (\ref{eq:many:costdiffo}) and (\ref{eq:many:costdiffm}) are meaningless. Observe that conditional on the label of the root $\phi$ the subtrees $\pwit_{\alpha,i},i\sim\phi$ are i.i.d., and the conditional distribution is same as the distribution of $\pwit_\alpha^m$ ($\pwit_\alpha^o$) given that $\phi$ has label $o$ ($m$). With this distribution structure of the subtrees, (\ref{eq:many:costdiffo}) and (\ref{eq:many:costdiffm}) motivate us to write a two-step {\em recursive distributional equation} (RDE): an equation in two distributions $\mu^o$ on $\Real$ and $\mu^m$ on $\RPlus$ that satisfy
\begin{equation}\label{eq:many:RDE}
\begin{split}
X^o&\eqdist\min_i\set{\alpha\xi_i-X^m_i},\\
X^m&\eqdist\min_i\set{\parl{\xi_i-X^o_i}^+},
\end{split}
\end{equation}
such that
\begin{enumerate}
\item $\set{\xi_i,i\ge 1}$ is a Poisson process of rate $1$ on $[0,\infty),$
\item $X^o$ has distribution $\mu^o,$
\item $X^m_i,i\ge 1$ are i.i.d.\ random variables with distribution $\mu^m,$ independent of the Poisson process $\set{\xi_i}$
\item $X^m$ has distribution $\mu^m,$
\item $X^o_i,i\ge 1$ are i.i.d.\ random variables with distribution $\mu^o,$ independent of the Poisson process $\set{\xi_i}$.
\end{enumerate}
We will call $(\mu^o,\mu^m)$ the solution to the RDE (\ref{eq:many:RDE}).

A recursive distributional equation is the primary tool that leads to the optimal solution on the PWIT $\pwit_\alpha$. The limit value of the expected optimal cost is obtained by calculations from the solution of the RDE. Recursive distributional equations come up in the related problems of random matching \cite{Ald2001} and edge-cover. Aldous and Bandyopadhyay \cite{AldBan2005} present a survey of this type of RDEs and introduce several RDE concepts that we will use later.

If $\mathcal{P}(S)$ denotes the space of probability measures on a space $S,$ a recursive distributional equation (RDE) is defined in \cite{AldBan2005} as a fixed-point equation on $\mathcal{P}(S)$ of the form
\begin{equation}\label{eq:many:rde}
X\eqdist g(\xi;(X_j,1\le j< N)),
\end{equation}
where $X_j,j\ge 1,$ are i.i.d.\ $S$-valued random variables having the same distribution as $X,$ and are independent of the pair $(\xi,N),$ $\xi$ is a random variable on some space, and $N$ is a random variable on $\Natural\cup\set{+\infty}$. $g$ is a given $S$-valued function. A solution to the RDE is a common distribution of $X, X_j,j\ge 1$ satisfying (\ref{eq:many:rde}).

The two-step equation (\ref{eq:many:RDE}) can easily be expressed in the form (\ref{eq:many:rde}) by composing the two steps of (\ref{eq:many:RDE}) into a single step. One can then solve for a fixed point distribution, say $\mu^o$, which then automatically yields $\mu^m$ via one of the steps of (\ref{eq:many:RDE}).

We can use the relation (\ref{eq:many:rde}) to construct a tree indexed stochastic process, say $X_{\bi i},{\bi i}\in\Vinf,$ which is called a recursive tree process (RTP) \cite{AldBan2005}. Associate to each vertex ${\bi i}\in\Vinf,$ an independent copy $(\xi_{\bi i},N_{\bi i})$ of the pair $(\xi,N),$ and require $X_{\bi i}$ to satisfy
\[
X_{\bi i}\eqdist g(\xi_{\bi i};(X_{{\bi i}.j},1\le j< N_{\bi i})),
\]
with $X_{\bi i}$ independent of $\set{(\xi_{\bi i^\prime},N_{\bi i^\prime}){\big\vert}\card{\bi i^\prime}<\card{\bi i}}$.
If $\mu\in\mathcal{P}(S)$ is a solution to the RDE~(\ref{eq:many:rde}), there exists a stationary RTP, i.e., each $X_{\bi i}$ is distributed as $\mu$. Such a process is called an invariant RTP with marginal distribution $\mu$.

\subsection{Solution to the RDE}
Write $F$ and $G$ for the complementary cdf of $\mu^o$ and $\mu^m$ respectively, i.e., $F(t)=\sP{\set{X^o>t}}$ and $G(t)=\sP{\set{X^m>t}}$. $F$ and $G$ satisfy the following:
\begin{gather}\label{eq:many:FfromG}
F(t)=\exp\left(-\int_{-t}^\infty \frac{G(z)}{\alpha}\rd z\right),\quad t\in\Real\\
\label{eq:many:GfromF}
G(t)=\begin{cases}\exp\left(-\int_{-t}^\infty F(z)\rd z\right)\quad\text{if } t\ge 0,\\
1\quad\text{if } t<0.
\end{cases}
\end{gather}

For $t\ge 0,$ we can write
\begin{align}
F(t)&=\exp\left(-\frac{1}{\alpha}\left(\int_{-t}^0\rd z - \int_0^\infty G(z)\rd z\right)\right)\nonumber\\
\label{eq:many:Ftwm}
&=e^{-w_m/\alpha}e^{-t/\alpha},
\end{align}
where $w_m=\sE{X^m}$.
Let $w_o=\sE{\parl{X^o}^+}$. By (\ref{eq:many:Ftwm}), we have
$w_o=\alpha e^{-w_m/\alpha}$. We can now rewrite (\ref{eq:many:Ftwm}) using $w_o$ as
\begin{equation}\label{eq:many:Ft}
F(t)=\frac{w_o}{\alpha}e^{-t/\alpha},\quad t\ge 0.
\end{equation}

From (\ref{eq:many:FfromG}) and (\ref{eq:many:GfromF}), $F$ and $G$ are differentiable for almost all $t\in\Real$. Differentiating (\ref{eq:many:FfromG}) and (\ref{eq:many:GfromF}), we get
\begin{align}\label{eq:many:Fprime}
F^\prime(-t)&=-\frac{1}{\alpha}F(-t)G(t),&\text{almost all }t\in\Real\\
\label{eq:many:Gprime}
G^\prime(t)&=-G(t)F(-t),&\text{almost all }t>0.
\end{align}
These two equations imply that
\[
\alpha F(-t)+G(t)=\text{constant},\quad t\ge 0.
\]
The validity of this equation can be extended to $t=0$ because $F$ is continuous for all $t$ and $G$ is right-continuous with left limits at $t=0$. Taking $t\To\infty$ in the above relation shows that the constant is equal to $\alpha,$ and so
\begin{equation}\label{eq:many:FG}
F(-t)=1-\frac{G(t)}{\alpha},\quad t\ge 0.
\end{equation}
Substituting this in (\ref{eq:many:Gprime}) results in the differential equation
\begin{equation}\label{eq:many:diffeq}
G^\prime(t)=-G(t)\left(1-\frac{G(t)}{\alpha}\right),\qquad\text{almost all }t>0.
\end{equation}
To get the boundary conditions take $t=0$ in (\ref{eq:many:GfromF}),(\ref{eq:many:Ft}), and (\ref{eq:many:FG}) to get
\begin{gather*}
G(0)=e^{-w_o},\\
F(0)=w_o/\alpha,\\
F(0)=1-G(0)/\alpha.
\end{gather*}
These imply $G(0)=\alpha-w_o$ and $\alpha=w_o+e^{-w_o}$. Since $\alpha>1,$ the latter equation determines a unique $w_o>0$.

Integration of the differential equation (\ref{eq:many:diffeq}) gives
\begin{equation*}
\begin{split}
G(t)&=\alpha \frac{G(0)e^{-t}}{\alpha-G(0)+G(0)e^{-t}}\\
&=\frac{\alpha}{1+w_oe^{w_o}e^t}
\end{split}
\end{equation*}
for $t\ge 0$. Using (\ref{eq:many:FG}), we get $F$ for negative values of $t$ as well.
We thus obtain the following expressions for the distributions $F$ and $G$:
\begin{align}\label{eq:many:distF}
F(t)&=\begin{cases}\frac{w_o}{\alpha}e^{-t/\alpha}&\text{if } t\ge 0,\\
1-\frac{1}{1+\gamma e^{-t}}&\text{if }t<0;
\end{cases}\\
G(t)&=\begin{cases}\frac{\alpha}{1+\gamma e^{t}}&\text{if }t\ge 0,\\ \label{eq:many:distG}
1&\text{if }t<0,
\end{cases}
\end{align}
where $\gamma=w_oe^{w_o}$.

For future use, we note that the density of $\mu^o$ is:
\begin{equation}\label{eq:many:fdensity}
 f(t)=\begin{cases}\frac{w_o}{\alpha^2}e^{-t/\alpha}&\text{if } t> 0,\\
\frac{\gamma e^{-t}}{(1+\gamma e^{-t})^2}&\text{if }t\le 0.
\end{cases}
\end{equation}

\section[Optimal involution invariant many-to-one matching on the infinite tree]{\texorpdfstring{Optimal involution invariant many-to-one matching on the infinite tree}{Optimal involution invariant many-to-one matching on the infinite tree}}\label{sec:many:theMatching}

The RDE (\ref{eq:many:RDE}) guides us to construct a many-to-one matching on $\pwit_\alpha$; the following Lemma forms the first step.

\begin{lemma}\label{lem:many:treeProcess}
There exists a process 
\[
\left(\pwit_\alpha,(\xi_e,e\in E(\pwit)),(\X{\overrightarrow{e}},\overrightarrow{e}\in \dirEdge{\pwit})\right)
\]
where $\pwit_\alpha$ is a modified PWIT with edge-lengths $\set{\xi_e,e\in E(\pwit)}$, and\newline $\set{\X{\overrightarrow{e}},\overrightarrow{e}\in \dirEdge{\pwit}}$ is a stochastic process satisfying the following properties.
\begin{enumerate}
    \item For each directed edge $(u,v)\in\dirEdge{\pwit},$
\begin{equation}\label{eq:many:Xrecursion}
\X{u,v}=\begin{cases}
\min\set{\len{v,w}-\X{v,w}:(v,w)\in\dirEdge{\pwit},w\neq u}\\
\hbox to 6.5cm{}\text{if }\mathsf{label}(v)=o,\\
\min\set{\diff{v}{w}:(v,w)\in\dirEdge{\pwit},w\neq u}\\
\hbox to 6.5cm{}\text{if }\mathsf{label}(v)=m.
\end{cases}
\end{equation}
\item If $(u,v)\in\dirEdge{\pwit}$ is directed away from the root of $\pwit,$ then $\X{u,v}$ has distribution $F$ conditioned on $v$ having label $o,$ and distribution $G$ conditioned on $v$ having label $m$; $F$ and $G$ are as in (\ref{eq:many:distF}) and (\ref{eq:many:distG}).
\item If $(u,v)\in\dirEdge{\pwit}$ the random variables $\X{u,v}$ and $\X{v,u}$ are independent.
\item For a fixed $z>0,$ conditional on the event that there exists an edge of length $z$ at the root, say $\set{\phi,v_z},$ and conditional on the root having label $m,$ the random variables $\X{\phi,v_z}$ and $\X{v_z,\phi}$ are independent random variables having distributions $F$ and $G$ respectively; if the root is conditioned to have label $o,$ the distributions are reversed.
\end{enumerate}
\end{lemma}
\begin{proof}
This Lemma is the analogue of Lemma~5.8 of \cite{AldSte2004} and Lemma~1 of \cite{KhaSun2012}. It has essentially the same proof, with appropriate changes to handle vertex labels. We omit the details.
\end{proof}

We use the process $\set{\X{\overrightarrow{e}}}$ to construct a many-to-one matching $\Mopt$ on $\pwit_\alpha$. For each vertex $v$ of $\pwit_\alpha,$ define a set 
\begin{equation} \label{eq:many:Moptdef}
\Mopt(v)=\begin{cases}
\argmin_{y\sim v}\set{\len{v,y}-\X{v,y}}&\text{if }\mathsf{label}(v)=o,\\
\argmin_{y\sim v}\set{\diff{v}{y}}&\text{if }\mathsf{label}(v)=m.
\end{cases}
\end{equation}
By continuity of the edge-weight distribution, when $v$ has label $o,$ almost surely $\Mopt(v)$ will be a singleton, whereas when $v$ has label $m,$ $\Mopt(v)$ may have more than one element -- e.g., all vertices $y$ for which the difference $\len{v,y}-\X{v,y}$ is negative.

Define
\[
\Mopt=\bigcup_v\set{\set{v,w}:w\in\Mopt(v)}.
\]
The following lemma proves that the collection $\Mopt$ yields a consistent many-to-one matching -- if $v$ with label $o$ picks $w$ with label $m$, then $w$ is the only $m$-labeled vertex that picks $v$.
\begin{lemma} \label{lem:many:consistentEC}
For any two vertices $v,w$ of $\pwit_\alpha,$ we have
\[
w\in\Mopt(v)\iff\len{v,w}<\X{v,w}+\X{w,v}.
\]
As a consequence,
\[
w\in\Mopt(v)\iff v\in\Mopt(w).
\]
\end{lemma}
\begin{proof}
Let $w\in\Mopt(v)$. Suppose $v$ has label $o$. 
Then 
\[
w=\argmin{\set{\len{v,y}-\X{v,y},y\sim v}}. 
\]
 This implies
\begin{align*}
\len{v,w}-\X{v,w}&=\min{\set{\len{v,y}-\X{v,y},y\sim v}}\\
&<\min{\set{\len{v,y}-\X{v,y},y\sim v,y\neq w}}=\X{w,v}.
\end{align*}
Conversely, if $\len{v,w}<\X{v,w}+\X{w,v}$ then using the definition of $\X{w,v},$ it follows that $w=\argmin{\set{\len{v,y}-\X{v,y},y\sim v}},$ and so $w\in\Mopt(v)$.

Now suppose $v$ has label $m$. Then there are two cases. (1) $\len{v,w}-\X{v,w}<0$. In this case $\len{v,w}<\X{v,w}+\X{w,v}$ since $\X{w,v}\ge 0$. Alternatively, (2) $\len{v,w}-\X{v,w}>0$, in which case, by (\ref{eq:many:Moptdef}), $w$ is the unique vertex where the minimum is attained. Hence
\[
0<\len{v,w}-\X{v,w}<\len{v,y}-\X{v,y}\quad\text{for all }y\sim v,y\neq w,
\]
and so
\[
\len{v,w}-\X{v,w}<\min{\set{\len{v,y}-\X{v,y},y\sim v,y\neq w}} =\X{w,v},
\]
which then yields $\len{v,w}<\X{v,w}+\X{w,v}$.
Again, it is easy to see that if $\len{v,w}<\X{v,w}+\X{w,v}$ then $w\in\argmin{\set{\diff{v}{y},y\sim v}},$ which implies $w\in\Mopt(v)$.

Thus, we have established the first statement of the lemma, which is
\[
w\in\Mopt(v)\iff\len{v,w}<\X{v,w}+\X{w,v}.
\]
The condition on the right-hand side above is symmetric in $v,w,$ and hence the second statement of the lemma is proved.
\end{proof}

\subsection{Evaluating the cost}
The next result is merely to demonstrate that the limiting expected minimum cost can be computed using numerical methods. 

Define 
\[
\calpha=-\polylog\parl{-\gamma^{-1}}-\frac{1}{2}\log^2\parl{1+\gamma^{-1}}+w_o\log\parl{1+\gamma^{-1}}+w_o,
\]
where $\polylog$ is the dilogarithm: $\polylog(z)=-\int_0^z\frac{\log(1-t)}{t}\rd t=\sum_{k=1}^\infty\frac{z^k}{k^2}=z \int_0^\infty t/(e^t-z) dt,$
$\gamma=w_o e^{w_o},$ and $w_o$ is the positive solution to $w_o+e^{-w_o}=\alpha$.
\begin{theorem}\label{thm:many:moptCost}
\begin{equation*}
\sEl{\sumver{v}{\Mopt}\len{\phi,v} \mid \mathsf{label}(\phi)=o } =\calpha.
\end{equation*}
\end{theorem}
\begin{proof}
Condition on the event that the root $\phi$ has label $o$. 
Now fix a $z>0,$ and condition on the event that there is a neighbor $v_z$ of $\phi$ with $\len{\phi,v_z}=z$. Call this event $E_z$. If we condition a Poisson process to have a point at some location, then the conditional process on removing this point is again a Poisson process with the same intensity. This shows that under $E_z,$ $\X{\phi,v_z}$ and $\X{v_z,\phi}$ have the distributions $G$ and $F$ respectively. They are also independent. Using these observations, the conditional expected cost can be written as
\begin{align*}
\lefteqn{\sEl{\sumver{v}{\Mopt}\len{\phi,v} \mid \mathsf{label}(\phi)=o }}&\\
&=\int_{z=0}^\infty z\sP{\set{z<\X{\phi,v_z}+\X{v_z,\phi}}} \frac{1}{\alpha} \rd z\\
&=\int_{z=0}^\infty\int_{x=-\infty}^\infty zG(z-x)f(x)\frac{1}{\alpha}\rd x\rd z\\
&=\int_{z=0}^\infty\int_{x=-\infty}^0 z\frac{1}{1+\gamma e^{z-x}}\frac{\gamma e^{-x}}{(1+\gamma e^{-x})^2}\rd x\rd z\\
&\qquad+\int_{z=0}^\infty\int_{x=0}^z z\frac{1}{1+\gamma e^{z-x}}\frac{w_o}{\alpha^2}e^{-x/\alpha}\rd x\rd z
+\int_{z=0}^\infty\int_{x=z}^\infty z f(x)\frac{1}{\alpha}\rd x\rd z\\
\intertext{(using the expression for the density of $\mu^o$ (\ref{eq:many:fdensity}))}
&=\int_{x=-\infty}^0 -\frac{\gamma e^{-x}}{(1+\gamma e^{-x})^2}\polylog(-\gamma^{-1}e^x)\rd x
+\frac{1}{\alpha^2}\int_{x=0}^\infty\int_{t=0}^\infty \frac{t+x}{1+\gamma e^t}w_oe^{-x/\alpha}\rd t\rd x\\
&\qquad+\frac{1}{\alpha}\int_{z=0}^\infty zF(z)\rd z \quad\parl{\text{taking }t=z-x}\\
&=\int_{s=\gamma }^\infty -\frac{1}{(1+s)^2}\polylog(-s^{-1})\rd s\\*
&\qquad+\frac{1}{\alpha^2}\int_{x=0}^\infty w_oe^{-x/\alpha}\parl{-\polylog(-\gamma^{-1})+x\log(1+\gamma^{-1})}\rd x\\*
&\qquad+\int_{z=0}^\infty z\frac{w_o}{\alpha^2}e^{-z/\alpha}\rd z \quad\parl{\text{taking }s=\gamma e^{-x}}\\
&=\frac{\polylog(-\gamma^{-1})}{1+\gamma }+\frac{1}{2}\log^2(1+\gamma^{-1})
-\frac{w_o}{\alpha} \polylog(-\gamma^{-1})\\
&\qquad+w_o\log(1+\gamma^{-1})+w_o\\
&=-\polylog(-\gamma^{-1})-\frac{1}{2}\log^2(1+\gamma^{-1})+w_o\log(1+\gamma^{-1})+w_o\\
&\quad\parl{\text{using }\frac{\alpha}{1+\gamma}=\alpha-w_o}.\qedhere
\end{align*}
\end{proof}

As $\alpha\To 1,$ the expected cost converges to $\pi^2/6,$ which is the limiting expected minimum cost of a standard matching on $\Knn$.

\subsection{Optimality in the class of involution invariant many-to-one matchings}
The local weak limit of a sequence of uniformly rooted random network satisfies a property called involution invariance \cite{AldSte2004, AldLyo2007}. Informally, involution invariance and its equivalent property called unimodularity specify that a distribution on the space of rooted networks $\graphs$
does not depend on the choice of the root.

Refer to Section~5 of \cite{KhaSun2012} for a discussion on involution invariance as it relates to our setting; a very general study of involution invariance is central to \cite{AldLyo2007}.

Let $\graphss$ be the space of isomorphism classes of connected geometric networks with a distinguished directed edge. An element of $\graphss$ $(G,o,x)$ will denote the isomorphism class corresponding to a network $G$ and distinguished vertex pair $(o,x)$, where $\set{o,x}$ is an edge in $G$. 
A probability measure $\mu$ on $\graphs$ is called involution invariant if it satisfies the following for all Borel $f:\graphss\To[0,\infty]$
\[
\int\sum_{x\in V(G)}f(G,o,x)\rd\mu(G,o)=\int\sum_{x\in V(G)}f(G,x,o)\rd\mu(G,o).
\]

We can represent a many-to-one matching $\matching$ on a graph $G$ by a membership map on the edge-set of $G$: $e\mapsto\ind_{\set{e\in\matching}}$. We call a random many-to-one matching on a random graph $G$ involution invariant if the distribution of $G$ with the above map on its edges is involution invariant. Following the argument in the last paragraph of Section~5 of \cite{KhaSun2012} we can restrict our attention to involution invariant many-to-one matchings on the PWIT. We now show that our candidate many-to-one matching $\Mopt$ has the minimum expected cost among such many-to-one matchings.

\begin{theorem}\label{thm:many:Mopt}
Let $\matching$ be an involution invariant many-to-one matching of the PWIT $\pwit_\alpha$. Write $\matching(\phi)$ for the set of vertices of $\pwit_\alpha$ adjacent to the root $\phi$ in $\matching$ (This set will be a singleton if label of $\phi$ is $o$). Then
\begin{equation*}
\sEl{\sumver{v}{\matching}\len{\phi,v}}\ge\sEl{\sumver{v}{\Mopt}\len{\phi,v}}.
\end{equation*}
\end{theorem}

Before we start to prove Theorem~\ref{thm:many:Mopt}, we examine a consequence of Theorem~\ref{thm:many:Mopt}.
\begin{corollary}\label{cor:many:Mopt}
\begin{equation*}
\sEl{\sumver{v}{\matching}\len{\phi,v} \mid \mathsf{label}(\phi)=o}\ge\sEl{\sumver{v}{\Mopt}\len{\phi,v} \mid \mathsf{label}(\phi)=o}.
\end{equation*}
\end{corollary}
\begin{proof}
  Involution invariance of $\matching$ implies that 
\[
\sEl{\sumver{v}{\matching}\len{\phi,v} \ind_{\set{\mathsf{label}(\phi)=o}} }=\sEl{\sumver{v}{\matching}\len{\phi,v} \ind_{\set{\mathsf{label}(\phi)=m}} }.
\]
Using that $\sP\set{\mathsf{label}(\phi)=o}=\frac{\alpha}{1+\alpha}$ and $\sP\set{\mathsf{label}(\phi)=m}=\frac{1}{1+\alpha},$ we get
\[
\sEl{\sumver{v}{\matching}\len{\phi,v}}=\frac{2\alpha}{1+\alpha}\sEl{\sumver{v}{\matching}\len{\phi,v} \mid \mathsf{label}(\phi)=o}.
\]
The same holds for $\Mopt$ too, and so the corollary follows from Theorem~\ref{thm:many:Mopt}.
\end{proof}

We will need some definitions for the proof of Theorem~\ref{thm:many:Mopt}.
For each directed edge $(v,w)$ of $\pwit_\alpha,$ if $v$ has label $o$ and $w$ has label $m$, define a random variable
\begin{equation} \label{eq:many:Ydef}
Y(v,w)=\min\set{\sum_{y\in A}(\len{w,y}-\X{w,y})\Big\vert\begin{matrix}A\subset\nbr{w}\setminus\set{v},\\A~\text{nonempty}\end{matrix}},
\end{equation}
where $\nbr{w}$ is the set of neighbors of $w$. It is easy to see that this random variable can be written as
\begin{equation*}
Y(v,w)=\begin{cases}
\mins{y}{w}{v}\set{\len{w,y}-\X{w,y}}\\
\qquad\qquad\qquad\text{if~}\len{w,y}-\X{w,y}\ge 0\text{~for~all~}\nbrneq{y}{w}{v}\\
\sum_{\nbrneq{y}{w}{v}}(\len{w,y}-\X{w,y})\ind_{\set{\len{w,y}-\X{w,y}<0}}\\
\qquad\qquad\qquad\text{otherwise.}
\end{cases}
\end{equation*} \label{eq:many:Ydefo}
Note that $\parl{Y(v,w)}^+=\X{v,w}$. If $v$ has label $o$ and $w$ has label $m$, define 
\begin{equation}
Y(w,v)=\X{w,v}.
\end{equation}

Suppose that $\sEl{\sumver{v}{\matching}\len{\phi,v}}<\infty$. Then $\matching(\phi)$ is a finite set with probability 1 because $\set{\len{\phi,v},v\sim\phi}$ are points of a Poisson process of rate 1 or $1/\alpha$. For such a many-to-one matching $\matching,$ define
\begin{equation}\label{eq:many:Adef}
A(\matching)=\sumver{v}{\matching}\X{\phi,v}+\maxver{v}{\matching}Y(v,\phi).
\end{equation}

The following two lemmas will be used to prove Theorem~\ref{thm:many:Mopt}.
\begin{lemma}\label{lem:many:D}
Let $\matching$ be a many-to-one matching on the PWIT such that 
\[
\sEl{\sumver{v}{\matching}\len{\phi,v}}<\infty.
\]
Then almost surely, 
\[
\sumver{v}{\matching}\len{\phi,v}\ge A(\matching).
\]
Furthermore,
\[
\sumver{v}{\Mopt}\len{\phi,v}=A(\Mopt).
\]
\end{lemma}

\begin{lemma}\label{lem:many:A}
Let $\matching$ be a many-to-one matching on the PWIT such that 
\[
\sEl{\sumver{v}{\matching}\len{\phi,v}} < \infty.
\] 
If $\matching$ is involution invariant, we have
$\sEl{A(\matching)}\ge \sEl{A(\Mopt)}$.
\end{lemma}

\begin{proof}[Proof of Theorem~\ref{thm:many:Mopt}]
If $\sEl{\sumver{v}{\matching}\len{\phi,v}}=\infty$ the statement of the theorem is trivially true. Assume that it is finite. We are now in a position to apply Lemmas \ref{lem:many:D} and \ref{lem:many:A} as follows to get the result:
\begin{align*}
\;\quad\quad\quad\quad\sEl{\sumver{v}{\matching}\len{\phi,v}}
&\ge \sEl{A(\matching)}\quad\text(Lemma~\ref{lem:many:D})\\
&\ge \sEl{A(\Mopt)}\quad\text(Lemma~\ref{lem:many:A})\\
&=\sEl{\sumver{v}{\Mopt}\len{\phi,v}}.\quad\text(Lemma~\ref{lem:many:D})\quad\quad\quad\quad\;\qedhere
\end{align*}
\end{proof}

Let us now complete the proofs of Lemmas \ref{lem:many:D} and \ref{lem:many:A}.
\begin{proof}[Proof of Lemma~\ref{lem:many:D}]
Suppose $\phi$ has label $o$. Then $\matching(\phi)$ is a singleton, say $\set{u}$. For $v\sim\phi,$ we have by (\ref{eq:many:Ydefo}) and (\ref{eq:many:Xrecursion}),
\[
Y(v,\phi)=\X{v,\phi}\le \len{\phi,y}-\X{\phi,y},y\sim\phi,y\neq v.
\]
If $v\neq u$, replacing $y$ with $u$ on the right-hand side, we get
\[
Y(v,\phi)\le\len{\phi,u}-\X{\phi,u}.
\]
Therefore,
\[
\maxver{v}{\matching}Y(v,\phi)\le\len{\phi,u}-\X{\phi,u}.
\]
This gives 
\[
\sumver{v}{\matching}\len{\phi,v}\ge A(\matching).
\]

If $\Mopt(\phi)=\set{u^*}$ then 
\[
u^*=\argmin_{y\sim\phi}\set{\len{\phi,y}-\X{\phi,y}}
\]
implies that, for any $v\neq u^*,v\sim\phi,$
\[
Y(v,\phi)=\X{v,\phi}=\len{\phi,u^*}-\X{\phi,u^*}.
\]
Therefore,
\[
\maxver{v}{\Mopt}Y(v,\phi)=\len{\phi,u^*}-\X{\phi,u^*},
\]
and so
\[
\sumver{v}{\Mopt}\len{\phi,v}=A(\Mopt).
\]

Now suppose $\phi$ has label $m$. Then, by (\ref{eq:many:Ydef}),
\[
Y(v,\phi)\le\sum_{y\in A}(\len{\phi,y}-\X{\phi,y})
\]
for all $A\subset\nbr{\phi}\setminus\set{v},$ $A$ nonempty.
In particular, for any $v\notin\matching(\phi),$ we can choose $A=\matching(\phi)$ to obtain
\[
Y(v,\phi)\le\sum_{y\in \matching(\phi)}(\len{\phi,y}-\X{\phi,y}).
\]
This implies
\begin{equation}\label{eq:many:maxY}
\maxver{v}{\matching} Y(v,\phi)\le\sum_{y\in \matching(\phi)}(\len{\phi,y}-\X{\phi,y}).
\end{equation}
Thanks to the finite expectation assumption in the lemma, $\matching(\phi)$ is a finite set almost surely, and so $\sum_{y \in \matching(\phi)} \X{\phi,y}$ is finite almost surely. Rearrangement of (\ref{eq:many:maxY}) then yields 
\[
\sumver{v}{\matching}\len{\phi,v}\ge A(\matching).
\]

We can write $\Mopt(\phi)$ as
\begin{equation} \label{eq:many:Moptalt}
\Mopt(\phi)=\argmin_A\set{\sum_{y\in A}(\len{\phi,y}-\X{\phi,y}):A\subset\nbr{\phi},A\text{ nonempty}}.
\end{equation}
From (\ref{eq:many:Ydef}) and (\ref{eq:many:Moptalt}), for any $v\notin\Mopt(\phi),$ we have
\begin{equation*}
Y(v,\phi)=\sum_{y\in \Mopt(\phi)}(\len{\phi,y}-\X{\phi,y}),
\end{equation*}
and hence
\[
\maxver{v}{\Mopt}Y(v,\phi)=\sum_{y\in\Mopt(\phi)}(\len{\phi,y}-\X{\phi,y}).
\]
It follows by rearrangement that
\[
\sumver{v}{\Mopt}\len{\phi,v}=A(\Mopt).\qedhere
\]
\end{proof}

\begin{proof}[Proof of Lemma~\ref{lem:many:A}]
Define
\begin{equation}\label{eq:many:Atilde}
\tilde{A}(\matching)=\sumver{v}{\matching}\X{v,\phi}+\maxver{v}{\matching}Y(v,\phi).
\end{equation}
$\tilde{A}(\matching)$ is similar to $A(\matching)$, but with the arguments of $X$ alone reversed.
We will prove Lemma~\ref{lem:many:A} by showing the following two results:
\begin{enumerate}
\item For an involution invariant many-to-one matching $\matching$
\begin{equation}\label{eq:many:lemA1}
\sEl{\tilde{A}(\matching)}=\sEl{A(\matching)}.
\end{equation}
\item Almost surely,
\begin{equation}\label{eq:many:lemA2}
\tilde{A}(\matching)\ge \tilde{A}(\Mopt).
\end{equation}
\end{enumerate}

We first prove (\ref{eq:many:lemA1}). First, by involution invariance of $\matching,$ we have
\begin{equation} \label{eq:many:Xsymm}
\sEl{\sumver{v}{\matching}\X{\phi,v}}=\sEl{\sumver{v}{\matching}\X{v,\phi}}.
\end{equation}
Indeed, the left-hand side equals
\[
\intG{\mu_\matching}{\X{\phi,v}}
\]
where $\mu_\matching$ is the probability measure on $\graphs$ corresponding to $\matching$.
By involution invariance, this equals 
\[
\intG{\mu_\matching}{\X{v,\phi}},
\]
which is equal to the right-hand side of (\ref{eq:many:Xsymm}).
Thanks to the finite expectation assumption of the lemma, we saw in the proof of Lemma~\ref{lem:many:D} that the term $\maxver{v}{\matching}Y(v,\phi)$ is finite almost surely. Now observe that $A(\matching)$ (respectively $\tilde{A}(\matching)$) is obtained by adding the almost surely finite random variable $\maxver{v}{\matching}Y(v,\phi)$ to the random variable which is the argument of the expectation in the left-hand side of (\ref{eq:many:Xsymm}) (respectively the right-hand side of (\ref{eq:many:Xsymm})). Taking expectation and using the equality in (\ref{eq:many:Xsymm}), we get (\ref{eq:many:lemA1}).

Now we will prove (\ref{eq:many:lemA2}). Suppose that $\phi$ has label $o$. Both $\matching(\phi)$ and $\Mopt(\phi)$ are singleton sets: suppose $\matching(\phi)=\set{u}$ and $\Mopt(\phi)=\set{u^*}$. If $u=u^*$ then $\tilde{A}(\matching)=\tilde{A}(\Mopt)$. Now suppose $u\neq u^*$. Then
\[
Y(v,\phi)=\X{v,\phi}=\len{\phi,u^*}-\X{\phi,u^*}
\]
for all $v\neq u^*,v\sim\phi$. For all such $v,$
\[
Y(v,\phi)=\X{u,\phi}\le\X{u^*,\phi}=Y(u^*,\phi).
\]
Consequently,
\begin{gather*}
\max_{v\neq u^*,v\sim\phi}Y(v,\phi)=\X{u,\phi},\\
\max_{v\neq u,v\sim\phi}Y(v,\phi)=\X{u^*,\phi}.
\end{gather*}
This implies that
\begin{gather*}
\tilde{A}(\Mopt)=\X{u^*,\phi}+\X{u,\phi},\\
\tilde{A}(\matching)=\X{u,\phi}+\X{u^*,\phi}.
\end{gather*}
Thus $\tilde{A}(\matching)=\tilde{A}(\Mopt)$.

Now suppose that $\phi$ has label $m$. First condition on the event 
\[
L_1=\set{\card{\Mopt(\phi)}>1}. 
\]
Observe that, under $L_1,$ $\len{\phi,y}-\X{\phi,y}<0,$ $y\sim\phi$ if and only if $y\in\Mopt(\phi),$ and there are at least two such $y$. Then, by (\ref{eq:many:Xrecursion}),
\begin{equation}\label{eq:many:Apart1}
\X{v,\phi}=0\;\text{for all }v\sim\phi.
\end{equation}
Also, from (\ref{eq:many:Ydef}) and (\ref{eq:many:Moptalt}), for any $v\sim\phi$, we have
\[
Y(v,\phi)\ge\sum_{y\in \Mopt(\phi)}(\len{\phi,y}-\X{\phi,y})=Y(w,\phi)
\]
if $w\notin\Mopt(\phi)$. This implies
\[
Y(v,\phi)\ge\maxver{w}{\Mopt}Y(w,\phi)\;\text{for all }v\sim\phi.
\]
In particular,
\begin{equation} \label{eq:many:Apart2}
\maxver{v}{\matching}Y(v,\phi)\ge \maxver{w}{\Mopt}Y(w,\phi).
\end{equation}
Combining (\ref{eq:many:Apart1}) and (\ref{eq:many:Apart2}) gives
\begin{equation*}
\sumver{v}{\matching}\X{v,\phi}+\maxver{v}{\matching}Y(v,\phi) \ge \sumver{v}{\Mopt}\X{v,\phi}+\maxver{v}{\Mopt}Y(v,\phi).
\end{equation*}
Thus $\tilde{A}(\matching)\ge\tilde{A}(\Mopt)$ under $L_1$.

Now consider the event $L_2=\set{\card{\Mopt(\phi)}=1}$. Let 
\begin{align*}
\Xone&=\min_{v\sim\phi}\parl{\len{\phi,v}-\X{\phi,v}}, \\
\Xtwo&=\mintwo{v\sim\phi}\parl{\len{\phi,v}-\X{\phi,v}},
\end{align*}
where $\mintwo{}$ stands for the second minimum.

Let $\Mopt(\phi)=\set{u^*}$. Then $\X{u^*,\phi}=\Xtwo$ (the second minimum is nonnegative), and for $v\in\matching(\phi)\setminus\Mopt(\phi),$ $\X{v,\phi}=\parl{\Xone}^+$. So we get
\begin{equation}\label{eq:many:Apart4}
\sumver{v}{\matching}\X{v,\phi}-\sumver{v}{\Mopt}\X{v,\phi}
=\sum_{v\in\matching(\phi)\setminus\Mopt(\phi)}\parl{\Xone}^+-\Xtwo\ind_{\set{u^*\notin\matching(\phi)}}.
\end{equation}
From (\ref{eq:many:Ydef}),
if $v\\neq u^*,$ then $Y(v,\phi)=\Xone$. Also $Y(u^*,\phi)=\Xtwo$. Since $\Xtwo\ge\Xone,$ we get
\[
\maxver{v}{\matching}Y(v,\phi)=\Xtwo\ind_{\set{u^*\notin\matching(\phi)}}+\Xone\ind_{\set{u^*\in\matching(\phi)}},
\]
and
\[
\maxver{v}{\Mopt}Y(v,\phi)=\Xone.
\]
Therefore,
\begin{equation}\label{eq:many:Apart5}
\maxver{v}{\matching}Y(v,\phi)-\maxver{v}{\Mopt}Y(v,\phi)=\left(\Xtwo-\Xone\right)\ind_{\set{u^*\notin\matching(\phi)}}.
\end{equation}
Adding (\ref{eq:many:Apart4}) and (\ref{eq:many:Apart5}), and canceling $\Xtwo\ind_{\set{u^*\notin\matching(\phi)}},$ we get
\begin{equation*}
\begin{split}
\sumver{v}{\matching}\X{v,\phi}+&\maxver{v}{\matching}Y(v,\phi)\\
&-\sumver{v}{\Mopt}\X{v,\phi}-\maxver{v}{\Mopt}Y(v,\phi)\\
&=\sum_{v\in\matching(\phi)\setminus\Mopt(\phi)}\parl{\Xone}^+-\Xone\ind_{\set{u^*\notin\matching(\phi)}}\\
&\ge 0,
\end{split}
\end{equation*}
where the last inequality follows because there exists a 
$v \in \matching(\phi) \setminus \Mopt(\phi)$ 
by virtue of our assumption that 
$\matching(\phi) \neq \Mopt(\phi)$.
Thus $\tilde{A}(\matching)\ge\tilde{A}(\Mopt)$ under $L_2$ as well.
\end{proof}

\section{Completing the lower bound}
\begin{theorem}\label{thm:many:lowerBd}
 Let $M_n^*$ be the optimal many-to-one matching on $\Knalpha$. Then
 \[
 \liminfn\frac{1}{n}\sEl{ \sum_{e\in M_n^*} \lennalpha{e}} \ge c_\alpha^*.
 \]
\end{theorem}
\begin{proof}
First observe that 
 \[
\sEl{ \sum_{e\in M_n^*} \lennalpha{e}} = n \sEl{ \sum_{\set{\phi,v}\in M_n^*} \len{\phi,v} \mid \mathsf{label}(\phi)=o},
 \]
because the root, conditioned to have label $o,$ is uniform among the $n$ vertices of the subset $A$.

The theorem can now be proved from Corollary~\ref{cor:many:Mopt} by following the exact steps of the proof of Theorem~9 in \cite{KhaSun2012}.
\end{proof}

\section{Endogeny} \label{sec:many:endogeny}
As we now move to establish the upper bound, our aim will be to show that BP run for large enough iterations  on $\Knalpha$ gives us a many-to-one matching that has cost close to the optimal for large $n$. We will show this by relating BP on $\Knalpha$ with BP on the PWIT through local weak convergence. It turns out that the $X$ process of Lemma~\ref{lem:many:treeProcess} generated from the RDE (\ref{eq:many:RDE}) arises as the limit of the message process on the PWIT (Theorem~\ref{thm:many:convOnT}). This result depends on an important property of the $X$ process called \emph{endogeny}, which was defined in \cite{AldBan2005}.
\begin{defnn}
An invariant RTP with marginal distribution $\mu$ is said to be endogenous if the root variable $X_\phi$ is almost surely measurable with respect to the $\sigma$-algebra $\sigma\parl{\set{\parl{\xi_{\bi i},N_{\bi i}} \mid {\bi i}\in\Vinf}}$.
\end{defnn}

Conceptually this means that the random variable at the root $X_\phi$ is a function of the ``data'' $\set{\parl{\xi_{\bi i},N_{\bi i}} \mid {\bi i}\in\Vinf}$, and there is no external randomness. This measurability also implies that if we restrict attention to a large neighborhood of the root, the root variable $X_\phi$ is almost deterministic given this neighborhood: the effect of the boundary vanishes. 

This sort of long range independence allows the use of belief propagation in conjunction with local weak convergence. When we run BP updates on a graph the number of iterations correspond to the size of the neighborhood around a vertex involved in the computation of the messages at that vertex. Endogeny gives us a way to say that this information is sufficient to make decisions at each vertex and jointly yield an almost optimal solution.

For a general RDE (\ref{eq:many:rde}), write $T:\mathcal{P}(S)\To\mathcal{P}(S)$ for the map induced by the function $g$ on the space of probability measures on $S$. We now define a bivariate map $T^{(2)}:\mathcal{P}(S\times S)\To\mathcal{P}(S\times S),$  which maps a distribution $\mu^{(2)}\in\mathcal{P}(S\times S)$ to the joint distribution of
\begin{equation*}
\begin{pmatrix}
g\parl{\xi;\parl{X^{(1)}_j,1\le j< N}}\\
g\parl{\xi;\parl{X^{(2)}_j,1\le j< N}}
\end{pmatrix}
\end{equation*}
where $\parl{X^{(1)}_j,X^{(2)}_j}_{j\ge 1}$ are independent with joint distribution $\mu^{(2)}$ on $S\times S,$ and the family of random variables $\parl{X^{(1)}_j,X^{(2)}_j}_{j\ge 1}$ are independent of the pair $(\xi,N)$.

It is easy to see that if $\mu$ is a fixed point of the RDE then the associated {\em diagonal measure} $\mudiag\coloneqq\text{Law}(X,X)$ where $X\sim\mu$ is a fixed point of the operator $T^{(2)}$.

Theorem~11c of \cite{AldBan2005}, reproduced below, provides a way of proving endogeny.
\begin{theorem}[Theorem 11c of \cite{AldBan2005}]\label{thm:many:endBVU}
For a Polish space $S,$ an invariant recursive tree process with marginal distribution $\mu$ is endogenous if and only if 
\[
  {T^{(2)}}^n(\mu\otimes\mu)\cvgdist\mudiag,                                                                                                      
\]
where $\mu\otimes\mu$ is the product measure.
\end{theorem}

We will use the above theorem to prove the following.
\begin{theorem}
 The invariant recursive tree process of Lemma~\ref{lem:many:treeProcess} arising from the solution to the RDE (\ref{eq:many:RDE}) is endogenous.
\end{theorem}
\begin{proof}
We first write an RDE for $X^o$. Let $g_1\parl{\parl{\xi_i}_{i\ge 1};\parl{x_i^o}_{i\ge 1}} = \min_i\parl{\xi_i-x_i^o}^+$.
Also define $g_2\parl{\parl{\xi_i}_{i\ge 1};\parl{x_i^m}_{i\ge 1}} = \min_i\parl{\xi_i-x_i^m}$.
Then define the map 
\[
g_o\parl{\parl{\parl{\xi_i}_{i\ge 1},\parl{\xi_{ij}}_{i,j\ge 1}};\parl{x_{ij}^o}_{i,j\ge 1}} 
= g_2\parl{\parl{\xi_i}_{i\ge 1};\parl{ g_1\parl{\parl{\xi_{ij}}_{j\ge 1};\parl{x_{ij}^o}_{j\ge 1}} }_{i\ge 1}}.
\]
We now have the RDE
\[
X^o\eqdist g_o\parl{\parl{\parl{\xi_i}_{i\ge 1},\parl{\xi_{ij}}_{i,j\ge 1}};\parl{X_{ij}^o}_{i,j\ge 1}},
\]
where 
\begin{enumerate}
 \item $\parl{\xi_i}_{i\ge 1}$ is a Poisson process of rate $1/\alpha$ on $\RPlus$,
 \item for each $i\ge 1,$ $\parl{\xi_{ij}}_{j\ge 1}$ is Poisson process of rate $1$ on $\RPlus$,
 \item ${X_{ij}^o, i,j\ge 1}$ are independent random variables having the same distribution as $X^o$, and
 \item the Poisson processes in (a) and (b) and the random variables $X_{i,j}$ in (c) are independent of each other.
\end{enumerate}
Call the induced mapping from $\mathcal{P}(\RPlus)$ to $\mathcal{P}(\RPlus)$ as $T$. The fixed point of this map is the distribution $\mu,$ with complementary cdf $F$ as in (\ref{eq:many:distF}).
We shall argue that 
\begin{equation} \label{eq:many:toshowendogeny}
  {T^{(2)}}^n(\mu\otimes\mu)\cvgdist\mudiag.
\end{equation}
Thus the invariant RTP associated with $g_o$ is endogenous. A similar argument shows that one can define an RDE $g_m$ and show that the corresponding invariant RTP associated with $g_m$ is endogenous. Recognizing that the layers of the two RTPs interlace, we easily see that the invariant RTP of Lemma~1 is endogenous. Let us now show (\ref{eq:many:toshowendogeny}).

Consider the bivariate distributional equations
\begin{align*}
\begin{pmatrix}
X^m_{(1)}\\
X^m_{(2)}
\end{pmatrix}
&\eqdist
\begin{pmatrix}
\min_i\parl{\xi_i-X_{i(1)}^o}^+\\
\min_i\parl{\xi_i-X_{i(2)}^o}^+\\
\end{pmatrix},\\
\begin{pmatrix}
X^o_{(1)}\\
X^o_{(2)}
\end{pmatrix}
&\eqdist
\begin{pmatrix}
\min_i\parl{\alpha\xi_i-X_{i(1)}^m}\\
\min_i\parl{\alpha\xi_i-X_{i(2)}^m}\\
\end{pmatrix},
\end{align*}
where $\set{\xi_i,i\ge 1}$ is a Poisson process of rate 1 on $\RPlus,$ $\set{\parl{X_{i(1)}^o,X_{i(2)}^o},i\ge 1}$ are i.i.d.\ random variables with joint complementary cdf $F_0^{(2)}$ independent of the Poisson process, 
$\set{\parl{X_{i(1)}^m,X_{i(2)}^m},i\ge 1}$ are i.i.d.\ random variables with joint complementary cdf $G_0^{(2)}$ independent of the Poisson process.
Let $\mathcal{P}^{(2)}$ denote the space of joint complementary cdfs of $\Real\times\Real$ valued random variables. Write $\Gamma^{(2)}:\mathcal{P}^{(2)}\To\mathcal{P}^{(2)}$ and $\varphi^{(2)}:\mathcal{P}^{(2)}\To\mathcal{P}^{(2)}$ for the maps defined by
\begin{align*}
\Gamma^{(2)}\parl{F_0^{(2)}}(x,y)&=\sP\set{X^m_{(1)}>x,X^m_{(2)}>y},\\
\varphi^{(2)}\parl{G_0^{(2)}}(x,y)&=\sP\set{X^o_{(1)}>x,X^o_{(2)}>y}.
\end{align*}

Let $F_0^{(2)}(x,y)=F(x)F(y)$ for all $x,y\in\Real,$ and define two sequences of joint distributions $\set{G_k^{(2)},k\ge 0}$ and $\set{F_k^{(2)},k\ge 0}$:
\begin{align*}
G_k^{(2)}&=\Gamma^{(2)}(F_k^{(2)}),\\
F_{k+1}^{(2)}&=\varphi^{(2)}(G_k^{(2)}).
\end{align*}
It is clear that the marginal distributions corresponding to $F_k^{(2)}$ and $G_k^{(2)}$ are $F$ and $G$ respectively for each $k$.
By Theorem~\ref{thm:many:endBVU}, the invariant RTP associated with $g_o$ is endogenous if $\set{F_k^{(2)},k\ge 1}$ converges to the joint distribution of the degenerate random variable that has both components equal. This is equivalent to the condition $F_k^{(2)}(x,x)\To F(x)$ as $k\To\infty$ for all $x\in\Real$.

From the bivariate RDE, if the joint distribution of $\parl{X^o_{i(1)},X^o_{i(2)}}$ is $F_k^{(2)}$ then for $x,y\ge 0,$
\begin{align*}
G_k^{(2)}(x,y)&=\sP\set{X^m_{(1)}>x,X^m_{(2)}>y}\\
&=\sP\left\{\text{No point of }\parl{\xi_i;X^o_{i(1)},X^o_{i(2)}}\right.\\
&\hbox to 1cm{}\left.\text{ in }\set{(z;u,v):z-u\le x\text{ or }z-v\le y}\right\}\\
&=\exp\parl{-\int_{z=0}^\infty\sP\set{z-X\le x\text{ or }z-Y\le y}\rd z}\\
&\hbox to 5cm{}\parl{(X,Y)\text{ has joint distribution }F_k^{(2)}}\\
&=\exp\parl{-\int_{z=0}^\infty(F(z-x)+F(z-y)-F_k^{(2)}(z-x,z-y))\rd z}\\
&=G(x)\exp\parl{-\int_{z=0}^\infty(F(z-y)-F_k^{(2)}(z-x,z-y))\rd z}.
\end{align*}
Write $g_k(x)=G_k^{(2)}(x,x)$ and $f_k(x)=F_k^{(2)}(x,x),x\in\Real$. Then
\begin{equation*}
g_k(x)=\begin{cases}G(x)\exp(-\int_{-x}^\infty(F(z)-f_k(z))\rd z),& x\ge 0,\\
                    1,& x<0.
       \end{cases}
\end{equation*}

Similar analysis with the equation for $\parl{X^o_{(1)},X^o_{(2)}}$ gives
\begin{equation}\label{eq:many:fk}
f_{k+1}(x)=F(x)\exp\parl{-\int_{-x}^\infty\frac{1}{\alpha}(G(z)-g_k(z))\rd z},x\in\Real.
\end{equation}

Define for $k\ge 0$,
\begin{align}
\beta_k(x)&=\int_{-x}^\infty(F(z)-f_k(z))\rd z,\label{eq:many:betak}\\
\gamma_k(x)&=\int_{-x}^\infty\frac{1}{\alpha}(G(z)-g_k(z)).\label{eq:many:gammak}
\end{align}
Then
\begin{align}
g_k(x)&=\begin{cases}G(x)\exp(-\beta_k(x)),& x\ge 0\\
                    1,& x<0,
       \end{cases}\label{eq:many:gkbetak}\\
f_{k+1}(x)&=F(x)\exp(-\gamma_k(x)),\;k\ge 0.\label{eq:many:fkgammak}
\end{align}

We now prove certain relations for the functions $\beta_k$ and $\gamma_k$.
\begin{lemma}
 For all $k\ge 0$:
 \begin{equation} \label{eq:many:betagammaprops}
\begin{split}
&\beta_k(x)\ge 0,\\
&\gamma_k(x)\ge 0,\\
&\beta_{k+1}(x)\le\beta_k(x),\text{ and}\\
&\gamma_{k+1}(x)\le\gamma_k(x) \text{ for all }x.
\end{split}
\end{equation}
\end{lemma}
\begin{proof}
 We will prove the Lemma by induction.
 
 $f_0(z)=F^2(z)\le F(z)$ for all $z$ implies that $\beta_0(x)\ge 0$ for all $x$. This implies $g_0(z)\le G(z)$ for all $z,$ and hence $\gamma_0(x)\ge 0$ for all $x$.
By (\ref{eq:many:fk}), we have 
\begin{gather*}
f_1(x)\ge F(x)\exp(-\int_{-x}^\infty\frac{1}{\alpha}G(z)\rd z)=f_0(x)\text{ for all }x,\\
\beta_1(x)\le\beta_0(x)\text{ for all }x\quad\text{by (\ref{eq:many:betak})},\\
g_1(x)\ge g_0(x)\text{ for all }x\quad\text{by (\ref{eq:many:gkbetak})},\\
\gamma_1(x)\le\gamma_0(x)\text{ for all }x\quad\text{by (\ref{eq:many:gammak})}.
\end{gather*}

We have showed that the equations (\ref{eq:many:betagammaprops}) hold for $k=0,$ now assume that they hold for some $k\ge 0$.

\begin{align*}
\gamma_k(x)\ge 0&\Rightarrow f_{k+1}(x)\le F(x)\text{ for all }x&&\text{by (\ref{eq:many:fkgammak})}\\
&\Rightarrow \beta_{k+1}(x)\ge 0\text{ for all }x&&\text{by (\ref{eq:many:betak})}\\
&\Rightarrow g_{k+1}(x)\le G(x)\text{ for all }x&&\text{by (\ref{eq:many:gkbetak})}\\
&\Rightarrow \gamma_{k+1}(x)\ge 0\text{ for all }x&&\text{by (\ref{eq:many:gammak})}.\\\\
\gamma_{k+1}(x)\le \gamma_k(x)&\Rightarrow f_{k+2}(x)\ge f_{k+1}(x)\text{ for all }x&&\text{by (\ref{eq:many:fkgammak})}\\
&\Rightarrow \beta_{k+2}(x)\le \beta_{k+1}(x)\text{ for all }x&&\text{by (\ref{eq:many:betak})}\\
&\Rightarrow g_{k+2}(x)\ge g_{k+1}(x)\text{ for all }x&&\text{by (\ref{eq:many:gkbetak})}\\
&\Rightarrow \gamma_{k+2}(x)\le\gamma_{k+1}(x)\text{ for all }x&&\text{by (\ref{eq:many:gammak})}.
\end{align*}
By induction, the properties (\ref{eq:many:betagammaprops}) hold for all $k\ge 0$.
\end{proof}

Because for each $x,$ the sequences $\set{\beta_k(x),k\ge 0}$ and $\set{\gamma_k(x),k\ge 0}$ are decreasing, the sequences of functions $\set{\beta_k, k\ge 0}$ and $\set{\gamma_k,k\ge 0}$ converge pointwise to functions, say $\beta$ and $\gamma$ respectively, as $k\To\infty$.
The requirement for endogeny that $f_k(x)\To F(x)$ for all $x$ as $k\To\infty$ reduces to showing that $\gamma(x)=0$ for all $x$. We will show that both $\beta$ and $\gamma$ are identically 0.

Using (\ref{eq:many:betak}) and (\ref{eq:many:fkgammak}), we have
\[
\beta_{k+1}(x)=\int_{-x}^\infty F(z)\parl{1-e^{-\gamma_k(z)}}\rd z.
\]
Since the term inside the integrand is bounded above by $F(z),$ and $\int_{-x}^\infty F(z)\rd z$ exists for all $x,$ by the dominated convergence theorem, we conclude that
\begin{equation}\label{eq:many:betagamma}
\beta(x)=\int_{-x}^\infty F(z)\parl{1-e^{-\gamma(z)}}\rd z.
\end{equation}
Similarly from (\ref{eq:many:gammak}) and (\ref{eq:many:gkbetak}), we have
\begin{equation}\label{eq:many:gammabeta}
\gamma(x)=\begin{cases}\int_{-x}^\infty \frac{G(z)}{\alpha}\parl{1-e^{-\beta(z)}}\rd z,&x<0\\
\int_{0}^\infty \frac{G(z)}{\alpha}\parl{1-e^{-\beta(z)}}\rd z,&x\ge 0
\end{cases}
\end{equation}

Differentiating (\ref{eq:many:betagamma}) and (\ref{eq:many:gammabeta}), we get (almost everywhere in the indicated sets)
\begin{align*}
\beta^\prime(x)&=F(-x)\parl{1-e^{-\gamma(-x)}},\;x\in\Real\\
\gamma^\prime(x)&=\frac{G(-x)}{\alpha}\parl{1-e^{-\beta(-x)}},\;x<0.
\end{align*}
Recall that for $x\ge 0$ (see (\ref{eq:many:Fprime}), (\ref{eq:many:Gprime}))
\begin{align*}
F^\prime(-x)&=-\frac{1}{\alpha}F(-x)G(x)\\
G^\prime(x)&=-G(x)F(-x).
\end{align*}

Define a function $\eta:\RPlus\To\RPlus$ as
\begin{equation}\label{eq:many:eta}
\eta(x)=\alpha F(-x)e^{\gamma(-x)}+G(x)e^{-\beta(x)}, x\ge 0.
\end{equation}
This function is continuous and differentiable on $(0,\infty)$.
Differentiating, we get (for $x>0$)
\begin{equation*}
\begin{split}
\eta^\prime(x)&=-\alpha F(-x)\gamma^\prime(-x)e^{\gamma(-x)}-\alpha F^\prime(-x)e^{\gamma(-x)}\\
&\quad -G(x)\beta^\prime(x)e^{-\beta(x)}+G^\prime(x)e^{-\beta(x)}\\
&=-F(-x)G(x)\parl{1-e^{-\beta(x)}}e^{\gamma(-x)}+F(-x)G(x)e^{\gamma(-x)}\\
&\quad-G(x)F(-x)\parl{1-e^{-\gamma(-x)}}e^{-\beta(x)}-G(x)F(-x)e^{-\beta(x)}\\
&=F(-x)G(x)e^{-\beta(x)}\parl{(e^{\gamma(-x)}+e^{-\gamma(-x)}-2}.
\end{split}
\end{equation*}
For any $y\in\Real,$ $e^y+e^{-y}\ge 2,$ hence
\[
\eta^\prime(x)\ge 0
\]
for all $x>0$.

Now observe that $\lim_{x\To\infty}\eta(x)=\alpha$.
By (\ref{eq:many:gammabeta}), $\gamma(x)=\gamma(0)$ for $x\ge 0$. Using this in (\ref{eq:many:betagamma}) gives
\[
\beta(0)=\int_0^\infty F(z)\parl{1-e^{-\gamma(0)}}\rd z=\parl{1-e^{-\gamma(0)}}w_o.
\]
Substituting values of $F(0)$ and $G(0)$ in (\ref{eq:many:eta}), and using the above relation between $\beta(0)$ and $\gamma(0),$ we get
\begin{equation*}
\begin{split}
\eta(0)&=w_oe^{\gamma(0)}+e^{-w_o}e^{-\beta(0)}\\
&=w_oe^{\gamma(0)}+e^{-w_o}e^{-\parl{1-e^{-\gamma(0)}}w_o}\\
&=w_oe^{\gamma(0)}+e^{-2w_o}e^{w_oe^{-\gamma(0)}}
\end{split}
\end{equation*}
Using $\gamma(0)\ge 0,$ it can be seen that $\eta(0)\ge w_o+e^{-w_o}=\alpha,$ with equality if and only if $\gamma(0)=0$.
We see that $\eta$ is a continuous increasing function on $[0,\infty)$ with $\eta(0)\ge \alpha$ and $\lim_{x\To\infty}\eta(x)=\alpha$. $\eta$ must then be a constant function with $\eta(0)=\alpha,$ and so $\gamma(0)=0$. For any $x\in\Real,$ $0\le\gamma(x)\le\gamma(0)$ implies that $\gamma$ must be identically 0. This proves that $f_k(x)$ converges to $F(x)$ for all $x\in\Real,$ and so the invariant RTP is endogenous.
\end{proof}

\section{Domain of attraction of the RDE map} \label{sec:many:domattr}
Consider the maps
\begin{align}
 \varphi G_0(t)&=\exp\left(-\int_{-t}^\infty\frac{G_0(z)}{\alpha}\rd z\right),\quad t\in\Real\label{eq:many:varphi}\\
 \Gamma F_0(t)&=\begin{cases}
 \exp\left(-\int_{-t}^\infty{F_0(z)}\rd z\right),\quad t\ge 0\label{eq:many:Gamma}\\
 1,\quad t<0.
               \end{cases}
\end{align}
Write $T=\Gamma\varphi$. We have $F=\varphi G,$ $G=\Gamma F,$ and $G=T G$ is the unique fixed point of $T$. 

We are interested in the domain of attraction of $G$. Specifically, we investigate the convergence of $T^k G_0$ for an arbitrary complementary cdf $G_0$. We will see that the domain of attraction is all those $G_0$ satisfying $\int_0^\infty G_0(z)\rd z<\infty$. We proceed via a sequence of Lemmas.

\begin{lemma} \label{lem:many:G0reg1}
Suppose $G_0$ is such that $\int_0^\infty G_0(t)\rd t<\infty$. Then,
\begin{enumerate}
 \item $TG_0(t)=1$ for $t<0,$
 \item $TG_0(0)<1,$
 \item $TG_0(t)$ is strictly decreasing for $t\ge 0,$
 \item $TG_0(\cdot)$ is continuous on $(0,\infty)$ (right-continuous at 0), and
 \item $\int_0^\infty TG_0(t)\rd t<\infty$.
\end{enumerate}
\end{lemma} 
\begin{proof}
 (a), (b), (c), and (d) follow directly from the definitions (\ref{eq:many:varphi}) and (\ref{eq:many:Gamma}).
 
 To get (e), observe that $F_0=\varphi G_0$ satisfies, for $z\le 0$, $F_0(z)\ge F_0(0) >0$ by the assumption that $\int_0^\infty G_0(t)\rd t<\infty$. Thus, for $t\ge 0$,
 \[
 \int_{-t}^\infty F_0(z)\rd z \ge tF_0(0).
 \]
 This implies, for $t\ge 0$,
 \[
 TG_0(t)=\exp\left(-\int_{-t}^\infty{F_0(z)}\rd z\right)\le e^{-tF_0(0)},
 \]
 which is integrable.
\end{proof}
Since, after one step of the iteration, $TG_0$ has properties (a), (b), (c), (d), and (e), we might assume, without loss of generality, that $G_0$ itself has these properties.

Define the transformation
\begin{equation}\label{eq:many:hatdef}
\widehat{G_0}(x)=x-\brackl{\log\parl{\frac{\alpha-G_0(x)}{\gamma G_0(x)}}}^+,\;x\ge 0.
\end{equation}
$\widehat{G_0}(x)$ is the largest shift at $x$ such that $G(x-\mathsf{shift})\le G_0(x)$. It has the following properties.
\begin{lemma}~\label{lem:many:hattrf}
\begin{enumerate}
\item $\widehat{G_0}(x)\le x$ with equality if and only if $G_0(x)\ge G(0)$.
\item $G_0(x)\ge G(x-\widehat{G_0}(x))$ with equality if and only if $G_0(x)\le G(0)$.
\end{enumerate}
\end{lemma}
\begin{proof}
 (a) Clearly $\widehat{G_0}(x)\le x$, because the second term that is subtracted in (\ref{eq:many:hatdef}) is $\ge 0$. Equality holds if and only if the logarithm in (\ref{eq:many:hatdef}) is $\le 0$
 \begin{align*}
  &\iff \frac{\gamma G_0(x)}{\alpha-G_0(x)} \ge 1\\
  &\iff G_0(x)\ge G(0)=\frac{\alpha}{1+\gamma }.
 \end{align*}
 
 (b) When $G_0(x)\ge G(0)$, we saw in (a) that $\widehat{G_0}(x)=x$. Hence in this case $G_0(x)\ge G(0)=G(x-\widehat{G_0}(x))$.
 We deal with the other case and the second part of (a) simultaneously by showing 
 \[
 G_0(x)\le G(0) \iff G_0(x) = G(x-\widehat{G_0}(x)).
 \]
This follows from the implications:
\begin{align*}
 G_0(x)\le G(0) &\iff \widehat{G_0}(x)=x-\log\parl{\frac{\alpha-G_0(x)}{\gamma G_0(x)}}\\
 &\iff x-\widehat{G_0}(x) = \log\parl{\frac{\alpha-G_0(x)}{\gamma G_0(x)}}\\
 &\iff G(x-\widehat{G_0}(x)) = \frac{\alpha}{1+\gamma \parl{\frac{\alpha-G_0(x)}{\gamma G_0(x)}} }=G_0(x). \qedhere
\end{align*}
\end{proof}

\begin{lemma}\label{lem:many:hatboundtoGineq}
Let $G_0$ satisfy the properties in Lemma~\ref{lem:many:G0reg1}. Let $x\ge 0$.
\begin{enumerate}
\item For $m\le 0,$ $\widehat{G_0}(x)\ge m$ if and only if $G_0(x)\ge G(x-m)$.
\item For $M\ge 0,$ $\widehat{G_0}(x)\le M$ if and only if $G_0(x)\le G(x-M)$.
\end{enumerate}
The above statements also hold when all the inequalities are replaced by strict inequalities.
\end{lemma}
\begin{proof}
Suppose that $\widehat{G_0}(x)\ge m$. Since $G$ is decreasing, we have $G(x-m)\le G(x-\widehat{G_0}(x))\le G_0(x)$ (the last inequality by Lemma~\ref{lem:many:hattrf}(b)).

Conversely, when $G_0(x)\ge G(x-m),$ either $G_0(x)=G(x-\widehat{G_0}(x)),$ which implies $\widehat{G_0}(x)\ge m,$ or by Lemma~\ref{lem:many:hattrf}(b) $G_0(x)>G(0),$ which by Lemma~\ref{lem:many:hattrf}(a) gives $\widehat{G_0}(x)=x\ge m$, because $x\ge 0$ and $m\le 0$.

Now let $\widehat{G_0}(x)\le M, M\ge 0$. Then for $x<M,$ $G_0(x)\le 1=G(x-M)$. For $x>M\ge \widehat{G_0}(x),$ by Lemma~\ref{lem:many:hattrf}(b), $G_0(x)<G(0),$ and so by Lemma~\ref{lem:many:hattrf}(a), $G_0(x)=G(x-\widehat{G_0}(x))\le G(x-M)$. The inequality extends to $x=M$ by right-continuity of $G_0$ and $G$ at $M$ and $0$ respectively.

Conversely, when $G_0(x)\le G(x-M),$ Lemma~\ref{lem:many:hattrf}(b) gives $G(x-\widehat{G_0}(x))\le G_0(x)\le G(x-M),$ and so $\widehat{G_0}(x)\le M$.

Also by strict monotonicity of $G_0(x)$ for $x\ge 0,$ all the inequalities can be made strict.
\end{proof}

\begin{lemma} \label{lem:many:Tmapboundshrink}
Suppose $G_0$ satisfies the properties in Lemma~\ref{lem:many:G0reg1}, and suppose further that there exist real numbers $m\le 0$ and $M\ge 0$ such that $m\le\widehat{G_0}(x)\le M$ for all $x\ge 0$. Then $m\le \widehat{TG_0}(x)\le M$ for all $x\ge 0$.
\end{lemma}
\begin{proof}
By Lemma~\ref{lem:many:hatboundtoGineq}, we have $G(z-m)\le G_0(z)\le G(z-M)$ for $z\ge 0$. But this implies $G(z-m)\le G_0(z)\le G(z-M)$ for all $z\in\Real$ as well, because for $z<0$, $G_0(z)=G(z-M)=1$.
Using these in (\ref{eq:many:varphi}), for all $t\in\Real,$ we have 
\begin{gather}
 \varphi G_0(t)\le \exp\left(-\int_{-t}^\infty\frac{G(z-m)}{\alpha}\rd z\right)=F(t+m),\label{eq:many:varphiineq1}\\
 \varphi G_0(t)\ge \exp\left(-\int_{-t}^\infty\frac{G(z-M)}{\alpha}\rd z\right)=F(t+M).\label{eq:many:varphiineq2}
\end{gather}
Then by (\ref{eq:many:Gamma}), for $x\ge 0$, we have
\begin{gather*}
 TG_0(x)=\Gamma\varphi G_0(x)\ge \exp\left(-\int_{-x}^\infty{F(t+m)}\rd t\right)=G(x-m),\\
 TG_0(x)=\Gamma\varphi G_0(x)\le \exp\left(-\int_{-x}^\infty{F(t+M)}\rd t\right)\le G(x-M),
\end{gather*}
where the last inequality follows from (\ref{eq:many:GfromF}) to handle the case $x<M$.

Lemma~\ref{lem:many:hatboundtoGineq} now completes the proof.
\end{proof}

The above Lemma shows that if $\widehat{G_0}$ is bounded then $\inf \widehat{T^kG_0}$ is increasing in $k$ and $\parl{\sup \widehat{T^kG_0}}^+$ is decreasing in $k,$ and so they converge to some $m^*$ and $M^*$ respectively ($m^*\le 0\le M^*$):
\begin{equation} \label{eq:many:Tmaplimits}
 \inf \widehat{T^kG_0}\uparrow m^* \text{ and } \parl{\sup \widehat{T^kG_0}}^+\downarrow M^* \text{ as } k\To\infty.
\end{equation}
We will show that $m^*=M^*=0$. (This will prove that both $\inf \widehat{T^kG_0}$ and $\sup \widehat{T^kG_0}$ converge to 0.)

First we show that the terms $\inf \widehat{T^kG_0}$ and $\parl{\sup \widehat{T^kG_0}}^+$ are strictly monotone in $k$ unless they are 0.
\begin{lemma}\label{lem:many:Tmapstrict}
Suppose $G_0$ satisfies the properties in Lemma~\ref{lem:many:G0reg1}. If $-\infty<\inf \widehat{G_0}<0$ then $\inf \widehat{TG_0}>\inf \widehat{G_0}$. If $0<\sup \widehat{G_0} <\infty$ then $\sup \widehat{TG_0} < \sup \widehat{G_0}$.
\end{lemma}
\begin{proof}
Let $M=\sup \widehat{G_0}>0$. Fix $a\in(0,M)$. By Lemma~\ref{lem:many:hatboundtoGineq}, $G_0(x)\le G(x-M)$. Using this in (\ref{eq:many:varphi}), we have, for $x\ge 0$
\begin{equation}\label{eq:many:varphistrictineq2}
 \begin{split}
 \varphi G_0(x)&\ge\exp\left(-\frac{1}{\alpha}\left(\int_{-x}^0G(z-M)\rd z+\int_{0}^aG_0(z)\rd z+\int_{a}^\infty G(z-M)\rd z\right)\right)\\
 &=\exp\left(-\int_{-x}^\infty\frac{G(z-M)}{\alpha}\rd z\right)\exp\left(\int_{0}^a\frac{G(z-M)-G_0(z)}{\alpha}\rd z\right)\\
 &=F(x+M)\exp\left(\int_{0}^a\frac{1-G_0(z)}{\alpha}\rd z\right)\\
 &=\kappa F(x+M),
 \end{split}
\end{equation}
where $\kappa>1$ (because $G_0(z)<1$ for $z\ge 0$).

Let $x\ge M$. Using (\ref{eq:many:varphiineq2}) to bound $\varphi G_0(z)$ for $z\in [-x,0)$ and (\ref{eq:many:varphistrictineq2}) to bound $\varphi G_0(z)$ for $z\ge 0$ in (\ref{eq:many:Gamma}), we get 
\begin{align}
 TG_0(x)&=\Gamma\varphi G_0(x)\notag\\
 &\le\exp\left(-\int_{-x}^0F(z+M)\rd z-\int_{0}^\infty\kappa F(z+M)\rd z\right)\notag\\
 &=\exp\left(-\int_{-x}^\infty{F(z+M)}\rd z\right)\exp\left(-(\kappa-1)\int_{0}^\infty F(z+M)\rd z\right)\notag\\
 &=\kappa^\prime G(x-M),\label{eq:many:Tstrictineq2}
\end{align}
where $\kappa^\prime<1$. Thus $TG_0(x)<G(x-M),$ and so by Lemma~\ref{lem:many:hatboundtoGineq}(b) for strict inequalities, $\widehat{TG_0}(x)<M$.

For $x\in[0,M),$ by Lemma~\ref{lem:many:hattrf}(a), we have $\widehat{TG_0}(x)\le x<M$. By continuity of $\widehat{TG_0}$, the supremum of $\widehat{TG_0}(x)$ over any compact subset of $[0,\infty)$ is strictly less than $M$. We now need to verify the strict inequality as $x\To\infty$. To show this, we have the following:
\begin{align*}
  \lim_{x\To\infty}\widehat{TG_0}(x)&=\lim_{x\To\infty} \left(x-\parl{\log\parl{\frac{\alpha-TG_0(x)}{\gamma TG_0(x)}}}^+\right) &&(\text{by definition (\ref{eq:many:hatdef})})\\
  &=\lim_{x\To\infty} \left(x-\log\parl{\frac{\alpha-TG_0(x)}{\gamma TG_0(x)}}\right) &&\\
  &&&\hskip-2cm\parl{\text{by }\lim_{x\To\infty}TG_0(x)=0<G(0)}\\
  &\le\lim_{x\To\infty}\left(x-\log\left(\frac{\alpha-\kappa^\prime G(x-M)}{\gamma \kappa^\prime G(x-M)}\right)\right) &&(\text{by (\ref{eq:many:Tstrictineq2})})\\
  &=\lim_{x\To\infty}\left(x-\log\left(\frac{\alpha-\kappa^\prime \frac{\alpha}{1+\gamma e^{x-M}}}{\gamma \kappa^\prime \frac{\alpha}{1+\gamma e^{x-M}}}\right)\right) &&\\
  &=\lim_{x\To\infty}\left(x-\log\left(\frac{1+\gamma e^{x-M}}{\gamma \kappa^\prime} - \frac{1}{\gamma}\right)\right) &&\\  
  &=M+\log\kappa^\prime<M.
\end{align*}

We use the same technique for the infimum. Suppose $m=\inf \widehat{G_0}<0$. Fix $b\in(m,0)$. Using $G_0(x)\ge G(x-m)$ for $x\ge 0$ (Lemma~\ref{lem:many:hatboundtoGineq}(a)) in (\ref{eq:many:varphi}), we have for $x\ge 0$
\begin{equation}\label{eq:many:varphistrictineq1}
 \begin{split}
 \varphi G_0(x)&\le\exp\left(-\frac{1}{\alpha}\left(\int_{-x}^b G(z-m)\rd z+\int_{b}^0 G_0(z)\rd z+\int_{0}^\infty G(z-m)\rd z\right)\right)\\
 &=\exp\left(-\int_{-x}^\infty\frac{G(z-m)}{\alpha}\rd z\right)\exp\left(\int_{b}^0\frac{G(z-m)-G_0(z)}{\alpha}\rd z\right)\\
 &=F(x+m)\exp\left(-\int_{b}^0\frac{1-G(z-m)}{\alpha}\rd z\right)\\
 &=\delta F(x+M),
 \end{split}
\end{equation}
where $\delta<1$. Using (\ref{eq:many:varphiineq1}) to bound $\varphi G_0(z)$ for $z<0$ and (\ref{eq:many:varphistrictineq1}) to bound $\varphi G_0(z)$ for $z\ge 0$ in (\ref{eq:many:Gamma}), we get for $x\ge 0$
\begin{equation*}
 \begin{split}
 TG_0(x)=\Gamma\varphi G_0(x)&\ge\exp\left(-\int_{-x}^0F(z+m)\rd z-\int_{0}^\infty\delta F(z+m)\rd z\right)\\
 &=\exp\left(-\int_{-x}^\infty{F(z+m)}\rd z\right)\exp\left((1-\delta)\int_{0}^\infty F(z+m)\rd z\right)\\
 &=\delta^\prime G(x-m),
 \end{split}
\end{equation*}
where $\delta^\prime>1$. This implies the strict inequality $\widehat{TG_0}(x)>m$ for all $x\ge 0$ (Lemma~\ref{lem:many:hatboundtoGineq}(a)). We now verify the inequality for $x\To\infty$, as in the case of supremum:
\begin{equation*}
 \begin{split}
\lim_{x\To\infty}\widehat{TG_0}(x)&=\lim_{x\To\infty}x-\parl{\log\parl{\frac{\alpha-TG_0(x)}{\gamma TG_0(x)}}}^+\\
  &\ge\lim_{x\To\infty}x-\log\parl{\frac{\alpha-\delta^\prime G(x-m)}{\gamma \delta^\prime G(x-m)}}\\
  &=m+\log\delta^\prime>m.
 \end{split}
\end{equation*}
\end{proof}

\begin{lemma} \label{lem:many:hatbounded}
Suppose $G_0$ satisfies the properties in Lemma~\ref{lem:many:G0reg1}. Then $\widehat{T^2G_0}$ is bounded.
\end{lemma}
\begin{proof}
Using $\int_0^\infty G_0(z)\rd z<\infty$ and $G_0(z)=1$ for $z<0$ in (\ref{eq:many:varphi}), we have
\begin{equation*}
 \varphi G_0(x)=\bigO(e^{-x/\alpha})\text{ as }x\To+\infty.
\end{equation*}
Then $\int_0^\infty\varphi G_0(z)\rd z<\infty,$ and by (\ref{eq:many:Gamma})
\begin{equation}\label{eq:many:order1}
 TG_0(x)=\bigTheta\left(e^{-\int_{-x}^{0}\varphi G_0(z)\rd z}\right)\text{ as }x\To+\infty.
\end{equation}
Since $\varphi G_0(z)$ is decreasing (and not identically 0), we have
\begin{equation*}
 TG_0(x)=\bigO(e^{-ux})\text{ as }x\To+\infty,
\end{equation*}
$u=\varphi G_0(0)>0$.
Using this in (\ref{eq:many:varphi}) gives
\begin{equation*}
 \varphi TG_0(x)=1-\bigO(e^{ux})\text{ as }x\To-\infty.
\end{equation*}
Now replacing $\varphi G_0$ by $\varphi TG_0$ in (\ref{eq:many:order1}), we get
\begin{equation*}
 T^2G_0(x)=\bigO(e^{-x})\text{ as }x\To+\infty.
\end{equation*}
Furthermore, using $\varphi G_0(z)\le 1$ for $z<0$ in (\ref{eq:many:order1}), we also get
\begin{equation*}
 T^2G_0(x)=\bigOmega(e^{-x})\text{ as }x\To+\infty,
\end{equation*}
and so 
\begin{equation}\label{eq:many:order2}
 T^2G_0(x)=\bigTheta(e^{-x})\text{ as }x\To+\infty.
\end{equation}

Let us now choose $K$ sufficiently large such that 
\begin{equation*}
 G(x+K)=\frac{\alpha}{1+\gamma e^{x+K}}\le T^2G_0(x) \le \frac{\alpha}{1+\gamma e^{x-K}}=G(x-K)
\end{equation*}
for all sufficiently large $x$. This, by Lemma~\ref{lem:many:hatboundtoGineq}, proves that $\widehat{T^2G_0}(x)$ is bounded.
\end{proof}
For future reference, let us also note that, by (\ref{eq:many:varphi}),
\begin{equation}\label{eq:many:order3}
 \varphi T^2G_0(x)=1-\bigOmega(e^{x})\text{ as }x\To-\infty.
\end{equation}

\begin{lemma} \label{lem:many:hatderivatives}
Suppose $\widehat{G_0}$ is bounded and $G_0$ satisfies the properties in Lemma~\ref{lem:many:G0reg1}. Then $\widehat{T^kG_0}$ is differentiable on $(0,\infty)$ except possibly at one point, and the derivatives are uniformly integrable for $k\ge 3$:
\begin{equation*}
 \sup_{k\ge 3}\int_{x>M}\lvert\widehat{T^kG_0}^\prime(x)\rvert\rd x\xrightarrow{M\To\infty}0.
\end{equation*}
\end{lemma}
\begin{proof}
 Let $x_0^k=\inf\set{x\ge 0:T^kG_0(x)\le G(0)}$. Then 
$\widehat{T^kG_0}(x)=x$ for $x\in[0,x_0^k],$ and $\widehat{T^kG_0}^\prime(x)=1,$ $x\in(0,x_0^k)$. The boundedness of $\widehat{G_0}(x)$ by $M$, say, implies by Lemma~\ref{lem:many:Tmapboundshrink} that $\widehat{T^kG_0}(x)$ are uniformly upper bounded by $M$ for all $x,$ and by Lemma~\ref{lem:many:hattrf}(a) $x_0^k\le M$. For $x>x_0^k$, by definition of $x_0^k$, we have $T^kG_0(x)\le G(0)$, in which case
\[
\widehat{T^kG_0}(x)=x-\log\parl{\frac{\alpha-G_0(x)}{\gamma G_0(x)}}.
\]
Consequently,
\begin{equation*}
 \widehat{T^kG_0}^\prime(x)=1+\frac{\alpha \parl{T^kG_0}^\prime(x)}{T^kG_0(x)\parl{\alpha-T^kG_0(x)}}.
\end{equation*}
By (\ref{eq:many:Gamma}),
\begin{equation*}
 \widehat{T^kG_0}^\prime(x)=1-\frac{\alpha \parl{\varphi T^{k-1}G_0}(-x)}{\alpha-T^kG_0(x)}.
\end{equation*}
As (\ref{eq:many:order2}) and (\ref{eq:many:order3}) hold respectively for $T^kG_0$ and $\varphi T^{k}G_0$ for all $k\ge 2$, the above equation gives
\begin{equation*}
 \widehat{T^kG_0}^\prime(x)=\bigO(e^{-x}).
\end{equation*}
This gives the uniform integrability of the derivatives.
\end{proof}

\begin{theorem}
For any complementary cdf $G_0$ satisfying $\int_0^\infty G_0(z)\rd z<\infty,$
\begin{equation*}
 T^kG_0 \To G \text{ pointwise as } {k\To\infty}.
\end{equation*}
The expectations also converge:
\begin{equation*}
\limk\int_0^\infty T^kG_0(x)\rd x=\int_0^\infty G(x)\rd x.
\end{equation*}
\end{theorem}
\begin{proof}
By Lemmas~\ref{lem:many:G0reg1}, \ref{lem:many:Tmapboundshrink}, \ref{lem:many:hatbounded}, and \ref{lem:many:hatderivatives}, for $k\ge 5,$ $\widehat{T^kG_0}$ are uniformly bounded and have derivatives that are uniformly integrable.

We prove that the limits in (\ref{eq:many:Tmaplimits}) satisfy $m^*=M^*=0$. The functions $T^kG_0,k\ge 1$ are bounded and 1-Lipshitz on $(0,\infty)$. By Arzela-Ascoli theorem this sequence is relatively compact with respect to compact convergence. There exists a subsequence that converges:
\begin{equation}\label{eq:many:Titercvg}
 T^{n_k}G_0\xrightarrow{k\To\infty}G_\infty.
\end{equation}
The uniform continuity of the function $y\mapsto \log(\frac{\alpha-y}{\gamma y})$ on every compact subset of $(0,1)$ implies that the transforms (which are uniformly bounded) converge (compact convergence):
\begin{equation*}
 \widehat{T^{n_k}G_0}\xrightarrow{k\To\infty}\widehat{G_\infty}.
\end{equation*}
The uniform integrability of the derivatives in Lemma~\ref{lem:many:hatderivatives} says that for any $\epsilon>0$, there exists $M$ sufficiently large such that $\widehat{T^{n_k}G_0}(x)\in (\widehat{T^{n_k}G_0}(M)-\epsilon,\widehat{T^{n_k}G_0}(M)+\epsilon)$ for all $x>M$ and for all $k\ge 3$. This extends the uniform convergence of $\widehat{T^{n_k}G_0}\To\widehat{G_\infty}$ in the compact set $[0,M]$ to $[0,\infty)$. The uniform convergence implies
\begin{gather*}
 \inf \widehat{G_\infty} =\limk\inf \widehat{T^{n_k}G_0} =m^*,\\
 \sup \parl{\widehat{G_\infty} }^+=\limk\parl{\sup \widehat{T^{n_k}G_0} }^+=M^*.
\end{gather*}
Using the dominated convergence theorem in (\ref{eq:many:varphi}) and (\ref{eq:many:Gamma}), the convergence (\ref{eq:many:Titercvg}) results in the compact convergence
\begin{equation*}
 T^{n_k+1}G_0\xrightarrow{k\To\infty}TG_\infty.
\end{equation*}
The same arguments now conclude that
\begin{gather*}
 \inf \widehat{TG_\infty }=m^*,\\
 \sup \parl{\widehat{TG_\infty }}^+=M^*.
\end{gather*}
Lemma~\ref{lem:many:Tmapstrict} allows this only if $m^*=M^*=0,$ and so $T^kG_0(x)\xrightarrow{k\To\infty}G(x)$ for all $x$.

Furthermore, there exists $M>0$ such that for all $k\ge 2,$ $\sup \widehat{T^kG_0}\le M$. This implies that $T^kG_0(x)\le G(x-M)$ for all $x\ge 0$ (Lemma~\ref{lem:many:hatboundtoGineq}), which by dominated convergence theorem gives convergence of the expectations:
\begin{equation*}
\limk\int_0^\infty T^kG_0(x)\rd x=\int_0^\infty G(x)\rd x.
\end{equation*}
\end{proof}

\section{Belief propagation} \label{sec:many:bp}
We now show that the belief propagation (BP) algorithm on $\Knalpha$ converges to a many-to-one matching solution that is asymptotically optimal as $n\To\infty$.

\subsection{Convergence of BP on the PWIT}
In this section we will prove that the messages on $\pwit$ converge, and relate the resulting matching with the matching $\Mopt$ of Section~\ref{sec:many:theMatching}.

The message process can essentially be written as
\begin{equation} \label{eq:many:BPonT}
\XT{\dot{v},v}{k+1}=\begin{cases}\min_{i\ge 1}\set{\lenT{v,v.i}-\XT{v,v.i}{k}},\quad\text{if label}(v)=o\\
                     \min_{i\ge 1}\set{\diffT{v}{v.i}{k}},\quad\text{if label}(v)=m.
                    \end{cases}
\end{equation}
where the initial messages $\XT{\dot{v},v}{0}$ are i.i.d.\ random variables (zero in the case of our algorithm; see (\ref{eq:many:bpInit})).

By the structure of $\pwit$ it is clear that -- the initial distribution being fixed -- the distribution of the messages $\XT{\dot{v},v}{k},v\in\Vinf$ depends only on $k$ and the label of the corresponding vertex $v$. Also, it can be seen from the analysis of RDE (\ref{eq:many:RDE}) in Section~\ref{sec:many:rde} that if we denote the complementary cdf of the  distribution of messages at some step $k^\prime$ by $F_0$ for vertices with label $o$ and $G_0$ for vertices with label $m,$ then after one update the respective complementary cdfs will be given by the maps (\ref{eq:many:varphi}) and (\ref{eq:many:Gamma}): $\varphi G_0$ and $\Gamma F_0$. Considering the messages $\XT{\dot{v},v}{k}$ only for vertices $v$ having label $m,$ the common distribution $G_0$ at step $k^\prime$ evolves to $TG_0$ at step $k^\prime+2$ ($T=\Gamma\varphi$); the distribution $\Gamma F_0$ at step $k^\prime+1$ changes to $T\Gamma F_0$ at $k^\prime+3$. The distribution sequence generated by the BP iterations (for $k\ge k^\prime$) is 
obtained by interleaving the two sequences generated by applying the map $T$ iteratively to $F_0$ and $\Gamma F_0$. By Section~\ref{sec:many:domattr}, the distributions along both the odd and even subsequences converge to $G,$ and consequently from the nature of the two step RDE, we conclude that the common distribution of $\XT{\dot{v},v}{k},$ $v$ having label $m$ converges to $G,$ and the common distribution of $\XT{\dot{v},v}{k},$ $v$ having label $o$ converges to $F$ as $k\To\infty$.

By the endogeny (Section~\ref{sec:many:endogeny}) of the RDE -- considering the evolution of messages at vertices of the same label -- the following result follows in the same way as Lemma~6 in \cite{KhaSun2012}.
\begin{lemma}\label{lem:many:endL2conv}
If the initial messages $\XT{\dot{v},v}{0},\mathsf{label}(v)=m$ are i.i.d.\ random variables with complementary cdf $G$ then the message process (\ref{eq:many:BPonT}) converges in $L^2$ to the process $X$ as $k\To\infty$.
\end{lemma}

We now prove that if the initial values are i.i.d.\ random variables with some arbitrary distribution (not necessarily one with complementary cdf $G$), then the message process (\ref{eq:many:BPonT}) does indeed converge to the unique stationary configuration. Of course, the initial condition of particular interest to us is the all zero initial condition (\ref{eq:many:bpInit}), but we will prove a more general result.

The following lemma will allow us to interchange limit and minimization while working with the updates on $\pwit$.

\begin{lemma} \label{lem:many:controlOnT}
Let $\set{\XT{\dot{v},v}{0}:\mathsf{label}(v)=o}$ be initialized to i.i.d.\ random variables with a distribution on $\Real$ having complementary cdf $F_0,$ and let\newline $\set{\XT{\dot{v},v}{0}:\mathsf{label}(v)=m}$ be initialized to i.i.d.\ random variables with a distribution on $\RPlus$ having complementary cdf $G_0$. Then for $v$ having label $o$, the map
\[
\mapT{(v)}{k}=\argmin_{u\sim v}\set{\lenT{v,u}-\XT{v,u}{k}}
\]
is a.s.\ well defined and singleton for all $k\ge 1,$ and
\[
\sup_{k\ge 1}\sP{\set{\argmin_{i\ge 1}\set{\lenT{v,v.i}-\XT{v,v.i}{k}}\ge i_0}}\To 0\quad\text{as }i_0\To\infty.
\]
Correspondingly, for $v$ having label $m$ the map
\[
\mapT{(v)}{k}=\argmin_{u\sim v}\set{\diffT{v}{u}{k}}
\]
is a.s.\ well defined and finite for all $k\ge 1,$ and
\[
\sup_{k\ge 1}\sP{\set{\max\argmin_{i\ge 1}\set{\diffT{v}{v.i}{k}}\ge i_0}}\To 0\quad\text{as }i_0\To\infty.
\]
\end{lemma}
The proof of this Lemma follows from the proofs of Lemma~7 of \cite{KhaSun2012} and Lemma~5.4 of \cite{SalSha2009}.

We are now in a position to prove the required convergence.

\begin{theorem} \label{thm:many:convOnT}
The recursive tree process defined by (\ref{eq:many:BPonT}) with i.i.d.\ initial messages (according to labels) converges to the unique stationary configuration in the following sense. For every $v\in\Vinf$ 
\[
\XT{v,v.i}{k}\cvgltwo \X{v,v.i}\quad\text{as }k\To\infty.
\]
Also, the decisions at the root converge, i.e.,  $\sP{\set{\mapT{(\phi)}{k}\neq\Mopt(\phi)}}\To 0$ as $k \To \infty$.
\end{theorem}
This can be proved by applying the proof of Theorem~13 of \cite{KhaSun2012} and Theorem~5.2 of \cite{SalSha2009} separately to the vertices having label $o$ and $m$.

\subsection{\texorpdfstring{Convergence of the update rule on $\Knalpha$ to the update rule on $\pwit$}{Convergence of the update rule on Kn,n/alpha to the update rule on T}}
We use from \cite{SalSha2009} the modified definition of local convergence applied to geometric networks with edge labels, i.e., networks in which each directed edge $(v,w)$ has a label $\lambda(v,w)$ taking values in some Polish space. For local convergence of a sequence of such labeled networks $G_1,G_2,\ldots$ to a labeled geometric network $G_\infty,$ we add the additional requirement that the rooted graph isomorphisms $\gamma_{n,\rho}$ satisfy
\[
\lim_{n\To\infty}\lambda_{G_n}(\gamma_{n,\rho}(v,w))=\lambda_{G_\infty}(v,w)
\]
for each directed edge $(v,w)$ in $\nbd{G_\infty}$.

Now we view the configuration of BP on a graph $G$ at the $k^{\text{th}}$ iteration as a labeled geometric network with the label on edge $(v,w)$ given by the pair 
\[
\left(\Xg{v,w}{G}{k},\ind_{\set{v\in\mapg{(w)}{G}{k}}}\right).
\]

With this definition, our convergence result can be written as the following theorem.

\begin{theorem} \label{thm:many:updateConv}
For every fixed $k\ge 0,$ the $k^\text{th}$ step configuration of BP on $\Knalpha$ converges in the local weak sense to the $k^\text{th}$ step configuration of BP on $\pwit$.
\begin{equation}\label{eq:many:updateConv}
\left(\Knalpha,\Xnalpha{v,w}{k},\ind_{\set{v\in\mapg{(w)}{\Knalpha}{k}}}\right) \cvglw\left(\pwit,\XT{v,w}{k},\ind_{\set{v\in\mapT{(w)}{k}}}\right)
\end{equation}
\end{theorem}
\begin{proof}
The proof of this theorem proceeds along the lines of the proof of Theorem~4.1 of \cite{SalSha2009}.

Consider an almost sure realization of the convergence $\Knalpha\To\pwit$. For such convergence to hold, the labels of the roots must match eventually, i.e., almost surely label of the root of $\Knalpha$ equals the label of the root of $\pwit$ for sufficiently large $n$.

Recall from Section~\ref{sec:many:lwConv} the enumeration of the vertices of $\pwit$ from the set $\Vinf$. We now recursively enumerate the vertices of $\Knalpha$ with multiple members of $\Vinf$. Denote the root by $\phi$. If $v\in\Vinf$ denotes a vertex of $\Knalpha,$ then $(v.1,v.2,\ldots,v.(h_v-1))$ denote the neighbors of $v$ in $\Knalpha$ ordered by increasing lengths of the corresponding edge with $v$; $h_v$ is the number of neighbors of $v$ ($\nbyalpha$ if $v$ has label $o$ and $n$ if $v$ has label $m$) . Then the convergence in (\ref{eq:many:updateConv}) is shown if we argue that
\begin{multline*}
\forall\set{v,w}\in\Einf\quad\Xnalpha{v,w}{k}\cvgpr\XT{v,w}{k}\quad\text{and}\quad\forall v\in\Vinf\quad\mapg{(v)}{\Knalpha}{k}\cvgpr\mapT{(v)}{k}\\
\text{as }n\To\infty. 
\end{multline*}

The above is trivially true for $k=0$. We write the update and decision rules, when $v$ has label $o,$ as
\begin{align*}
\Xnalpha{w,v}{k+1}&=\min_{u\in\set{v.1,\ldots,v.(\nbyalpha-1),\dot{v}}\setminus\set{w}} \set{\lennalpha{v,u}-\Xnalpha{v,u}{k}}\quad\text{and}\\
\mapg{(v)}{\Knalpha}{k}&=\argmin_{u\in\set{v.1,\ldots,v.(\nbyalpha-1),\dot{v}}}\set{\lennalpha{v,u}-\Xnalpha{v,u}{k}},
\end{align*}
and when $v$ has label $m,$ as
\begin{align*}
\Xnalpha{w,v}{k+1}&=\min_{u\in\set{v.1,\ldots,v.(n-1),\dot{v}}\setminus\set{w}} \set{\parl{\lennalpha{v,u}-\Xnalpha{v,u}{k}}^+}\quad\text{and}\\
\mapg{(v)}{\Knalpha}{k}&=\argmin_{u\in\set{v.1,\ldots,v.(n-1),\dot{v}}}\set{\parl{\lennalpha{v,u}-\Xnalpha{v,u}{k}}^+},
\end{align*}

We want to use the convergence of each term on the right-hand side inductively to conclude the convergence of the term on the left. The minimum taken over an unbounded number of terms as $n\To\infty$ creates a problem. However the following lemma allows us to restrict attention to a uniformly bounded number of terms for each $n$ with probability as high as desired, and hence obtain convergence in probability for each $k\ge 0$.
\end{proof}

\begin{lemma} \label{lem:many:unifControl}
For all $k\ge 0,$ when $v\in\Vinf$ has label $o$
\[
\lim_{i_0\To\infty}\limsupn\sP{\set{\argmin_{1\le i\le\nbyalpha-1}\set{\lennalpha{v,v.i}-\Xnalpha{v,v.i}{k}}\ge i_0}}=0,
\]
and when $v$ has label $m$
\[
\lim_{i_0\To\infty}\limsupn\sP{\set{\max\argmin_{1\le i\le n-1}\set{\parl{\lennalpha{v,v.i}-\Xnalpha{v,v.i}{k}}^+}\ge i_0}}=0.
\]
\end{lemma}
\begin{proof}
The proof is the same as the proof of Lemma~4.1 of \cite{SalSha2009}. Also see Lemma~8 in \cite{KhaSun2012}.
\end{proof}

\subsection{Completing the upper bound - Proof of Theorem~\ref{thm:many:bpconv}}

By Theorem~\ref{thm:many:convOnT}, $\mapT{(\phi)}{k}\cvgpr\Mopt(\phi)\text{ as }k\To\infty$. It follows that 
\begin{equation}\label{eq:many:costconvk}
\sum_{v\in\mapT{(\phi)}{k}}\lenT{\phi,v}
\cvgpr\sum_{v\in\Mopt{(\phi)}}\lenT{\phi,v}\quad\text{as }k\To\infty.
\end{equation}
We now prove convergence in expectation.
Observe that
\[
v\in\mapT{(\phi)}{k} \Rightarrow \lenT{\phi,v}-\XT{\phi,v}{k}\le\diffT{\phi}{1}{k}\le\lenT{\phi,1}.
\]
By (\ref{eq:many:BPonT}), $\XT{\phi,v}{k}\le\lenT{v,v.1}$. Thus,
\begin{equation}\label{eq:many:dominatek}
v\in\mapT{(\phi)}{k} \Rightarrow \lenT{\phi,v}\le\lenT{\phi,1}+\lenT{v,v.1}.
\end{equation}
This implies
\begin{equation*}
\sum_{v\in\mapT{(\phi)}{k}}\lenT{\phi,v}\le \lenT{\phi,1}+\sum_{i\ge 2}\lenT{\phi,i}\ind_{\set{\lenT{\phi,i}\le\lenT{\phi,1}+\lenT{i,i.1}}}.
\end{equation*}
It can be verified that the sum on the right-hand side in the above equation is an integrable random variable. Equation~(\ref{eq:many:costconvk}) and the dominated convergence theorem give
\begin{equation*}
\limk\sEl{\sum_{v\in\mapT{(\phi)}{k}}\lenT{\phi,v}}
=\sEl{\sum_{v\in\Mopt{(\phi)}}\lenT{\phi,v}}.
\end{equation*}
Involution invariance gives (see Corollary~\ref{cor:many:Mopt})
\begin{align}
\limk\sEl{\sum_{v\in\mapT{(\phi)}{k}}\lenT{\phi,v}\mid\mathsf{label}(\phi)=o}
&=\sEl{\sum_{v\in\Mopt{(\phi)}}\lenT{\phi,v}\mid\mathsf{label}(\phi)=o}\notag\\
&=c^*_\alpha.\label{eq:many:expCostLimitk}
\end{align}
where the last equality follows from Theorem~\ref{thm:many:moptCost}.

By Theorem~\ref{thm:many:updateConv} and Lemma~\ref{lem:many:unifControl}, using the definition of local weak convergence, we have
\begin{equation}\label{eq:many:dominaten}
\sum_{v\in\mapnalpha{(\phi)}{k}}\lennalpha{\phi,v}
\cvgpr\sum_{v\in\mapT{(\phi)}{k}}\lenT{\phi,v}\quad\text{as }n\To\infty.
\end{equation}
We now apply the arguments that lead to (\ref{eq:many:dominatek}) to $\mapnalpha{(\phi)}{k},$ and obtain
\begin{equation*}
v\in\mapnalpha{(\phi)}{k} \Rightarrow \lennalpha{\phi,v}\le\lennalpha{\phi,1}+\lennalpha{v,v.1}.
\end{equation*}
For any two vertices $u,v$ of $\Knalpha,$ define $S_n(u,v)=\min_{w\neq u,v}\lennalpha{u,w}$. Then for a vertex $v$ of $\Knalpha,$ $\lennalpha{\phi,1}\le S_n(\phi,v)$ and $\lennalpha{v,v.1}\le S_n(v,\phi)$. This gives
\begin{equation*}
v\in\mapnalpha{(\phi)}{k} \Rightarrow \lennalpha{\phi,v}\le S_n(\phi,v)+S_n(v,\phi).
\end{equation*}
Consequently,
\begin{equation}\label{eq:many:dominaten1}
\sum_{v\in\mapnalpha{(\phi)}{k}}\lennalpha{\phi,v}\le \sum_v\lennalpha{\phi,v}\ind_{\set{\lennalpha{\phi,v}\le S_n(\phi,v)+S_n(v,\phi)}}.
\end{equation}

Condition that $\phi$ has label $o$. Observe that $\lennalpha{\phi,v},S_n(\phi,v),\text{ and }S_n(v,\phi)$ are independent exponential random variables with means $n,n/(\nbyalpha-1),$ and $n/(n-1)$ respectively. So we can write
\begin{equation*}
\begin{split}
\lefteqn{\sEl{\lennalpha{\phi,v}\ind_{\set{\lennalpha{\phi,v}\le S_n(\phi,v)+S_n(v,\phi)}} \mid \mathsf{label}(\phi)=o}}\qquad&\\
&=\int_0^\infty\int_0^x\frac{t}{n}e^{-t/n}\rd t\left(\frac{(\nbyalpha-1)(n-1)}{n(n-\nbyalpha)}\right) \left(e^{-\left(\frac{\nbyalpha-1}{n}\right)x}-e^{-\left(\frac{n-1}{n}\right)x}\right)\rd x\\
&=\left(\frac{n}{\nbyalpha}+1\right)\left(\frac{1}{\nbyalpha}+\frac{1}{n}-\frac{1}{n\nbyalpha}\right)-\frac{1}{\nbyalpha}.
\end{split}
\end{equation*}
Summing over all neighbors of $\phi,$ we get
\begin{multline}\label{eq:many:dominaten2}
\sEl{\sum_v\lennalpha{\phi,v}\ind_{\set{\lennalpha{\phi,v}\le S_n(\phi,v)+S_n(v,\phi)}}}\\
=\left(\frac{n}{\nbyalpha}+1\right)\left(1+\frac{\nbyalpha}{n}-\frac{1}{n}\right)-1,
\end{multline}
which converges to $1+\alpha+1/\alpha$ as $n\To\infty$.

Using local weak convergence, we can see that
\begin{multline*}
\sum_v\lennalpha{\phi,v}\ind_{\set{\lennalpha{\phi,v}\le S_n(\phi,v)+S_n(v,\phi)}}\\
\cvgpr\lenT{\phi,1}+\sum_{i\ge 2}\lenT{\phi,i}\ind_{\set{\lenT{\phi,i}\le\lenT{\phi,1}+\lenT{i,i.1}}}.
\end{multline*}
It can be verified that the expectation of the random variable on the right-hand side above equals $1+\alpha+1/\alpha$. Using this with (\ref{eq:many:dominaten}), (\ref{eq:many:dominaten1}), and (\ref{eq:many:dominaten2}), the generalized dominated convergence theorem yields
\begin{multline}\label{eq:many:expCostLimitn}
\limn\sEl{\sum_{v\in\mapnalpha{(\phi)}{k}}\lennalpha{\phi,v}\mid\mathsf{label}(\phi)=o}\\
=\sEl{\sum_{v\in\mapT{(\phi)}{k}}\lenT{\phi,v}\mid\mathsf{label}(\phi)=o}.
\end{multline}
Combining (\ref{eq:many:expCostLimitn}) and (\ref{eq:many:expCostLimitk}) gives
\begin{equation}\label{eq:many:expCostLim}
\limk\limn\sEl{\sum_{v\in\mapnalpha{(\phi)}{k}}\lennalpha{\phi,v}\mid\mathsf{label}(\phi)=o}= c_\alpha^*.
\end{equation}
Use of
\[
\sEl{\sum_{e\in\matching(\mapnalpha{}{k})}\lennalpha{e}}=n \sEl{\sum_{v\in\mapnalpha{(\phi)}{k}}\lennalpha{\phi,v}\mid\mathsf{label}(\phi)=o}
\]
with (\ref{eq:many:expCostLim}) completes the proof of Theorem~\ref{thm:many:bpconv}.

\subsection{Completing the proof of Theorem~\ref{thm:many:limit}}
As remarked earlier, $\matching(\mapnalpha{}{k}),$ unlike its limit $\Mopt,$ is not a many-to-one matching. The matching condition is violated for a vertex $v\in A$ if there are more than one vertices $w$ such that $v\in\mapnalpha{(w)}{k}$. We now show how to construct a valid many-to-one matching from $\matching\bigl(\mapnalpha{}{k}\bigr)$ with little increase in cost.

\begin{theorem} \label{thm:many:patchingmatching}
 Fix $\epsilon>0$. There exist many-to-one matchings $\matching_n$ on $\Knalpha$ such that for all sufficiently large $k$
 \[
\limsupn n^{-1} \sEl{\sum_{e\in\matching_n}\lennalpha{e}} \le \limn n^{-1} \sEl{\sum_{e\in\matching\bigl(\mapnalpha{}{k}\bigr)}\lennalpha{e}} +\epsilon.
\]
\end{theorem}
\begin{proof}
 Define a partial matching $\tilde{\matching}\bigl(\mapnalpha{}{k}\bigr)=\set{\set{v,w}\mid w\in\mapnalpha{(v)}{k}, v\in A}$. $\tilde{\matching}\bigl(\mapnalpha{}{k}\bigr)$ is a subset of $\matching\bigl(\mapnalpha{}{k}\bigr),$ therefore, its cost denoted $\mathsf{cost}_n^k$ is bounded above by the cost of $\matching\bigl(\mapnalpha{}{k}\bigr)$. Every vertex in $A$ has exactly one edge in $\tilde{\matching}\bigl(\mapnalpha{}{k}\bigr),$ but some vertices in $B$ may be unmatched. Call this set of unmatched vertices $U_n^k$.
 By Theorem~\ref{thm:many:updateConv}, for a vertex $w\in B,$ 
 \[
  \sP\set{w\in U_n^k}\To\sP\set{\phi\notin\bigcup_{v\sim\phi}\mapT{(v)}{k}\mid\mathsf{label}(\phi)=m}\eqqcolon r
 \]
as $n\To\infty$. We have $r\To 0$ as $k\To\infty$ (Theorem~\ref{thm:many:convOnT}).
For any $\delta>0,$ we then have $\sE\card{U_n^k}\le (r+\delta)n/\alpha$ for all sufficiently large $n$.
 
 Section~3.3 of \cite{Was2012} describes a method (which is based on \cite{Fri2004}) for constructing a perfect matching in the random complete graph $K_n$ from a {\em diluted }matching that has a small fraction of vertices unmatched, with a small increase in cost. Our case is much simpler because more than one edge is allowed for a vertex in $B$. We start with the first step of the proof in \cite{Was2012}, and deviate afterwards.
 
 In $\Knalpha$ substitute the single edge between every pair of vertices with a countably infinite number of edges and take their weights according to a Poisson process on $\RPlus$ with rate $1/n$. The weight of the smallest edge between every pair of vertices is still distributed as the original (exponential with mean $n$), and so the optimal many-to-one matching in this new network is same as that in $\Knalpha$. Independently color every edge red with probability $1-p$ and green with probability $p$ (we will later take $p$ to 0). Obtain a partial matching $\tilde{\matching}\bigl(\mapnalpha{}{k}\bigr)$ on this network using only red edges. The restriction to red edges just affects the scaling of the cost: the weight of the smallest red edge between a vertex pair has exponential distribution with mean $n/(1-p),$ and so the cost of this partial matching is $\mathsf{cost}_n^k/(1-p)$.
 
For every matched vertex $w$ in $B\setminus U_n^k$ leave out one vertex in $\mapnalpha{(w)}{k}$ (choice made arbitrarily if size more than 1) and collect the remaining vertices in a set $S_n^k$. There are exactly $n-(\nbyalpha-\card{U_n^k})$ vertices in this set. We can complete the many-to-one matching by greedily matching every vertex in $U_n^k$ with a distinct vertex in $S_n^k$ and dropping the original matched edge of that vertex. We will do this by using only green edges. The conditional expected cost of these green edges, given the coloring and the weight of the red edges, is at most $\frac{\card{U_n^k}}{\card{S_n^k}-\card{U_n^k}+1}\frac{n}{p}\le \frac{\card{U_n^k}}{(1-1/\alpha)p}$. We get this because the choice of the sets $U_n^k$ and $S_n^k$ is governed by the red edges, whose weights are independent of those of the green edges.

The expected cost of the constructed many-to-one matching is at most
\begin{align*}
 \sEl{\frac{\mathsf{cost}_n^k}{(1-p)}+\frac{\card{U_n^k}}{(1-1/\alpha)p}}
 \le \frac{\sEl{\mathsf{cost}_n^k}}{(1-p)}+\frac{(r+\delta)n}{(\alpha-1)p}
\end{align*}
for large $n$.

Since $r\To 0$ as $k\To\infty$ and $\delta>0,$ $p>0$ are arbitrary, we can bound the expected cost by $\sEl{\mathsf{cost}_n^k}+\epsilon n$ for large $k$ and large $n$.
\end{proof}

By Theorem~\ref{thm:many:bpconv} for any $\epsilon>0,$ we can find $k$ large such that 
\[
\limn n^{-1} \sEl{\sum_{e\in\matching(\mapnalpha{}{k})}\lennalpha{e}} \le c_\alpha^*+\epsilon.
\]
This combined with Theorem~\ref{thm:many:patchingmatching} gives
\[
\limsupn n^{-1} \sE{\Mn}\le c_\alpha^*+2\epsilon.
\]
Since $\epsilon$ is arbitrary, we get the upper bound
\[
\limsupn n^{-1} \sE{\Mn}\le c_\alpha^*.
\]
We have the lower bound from Theorem~\ref{thm:many:lowerBd}.
This completes the proof of Theorem~\ref{thm:many:limit}.\qed


Observe that for any $\epsilon>0,$ there exist $K_\epsilon$ and $N_\epsilon$ such that for all $k\ge K_\epsilon$ and $n\ge N_\epsilon,$ we have
\[
 n^{-1}\sEl{\sum_{e\in\matching(\mapnalpha{}{k})}\lennalpha{e}}\le c_\alpha^*+\epsilon.
\]
Thus for large $n$ the BP algorithm gives a solution with cost within $\epsilon$ of the optimal value in $K_\epsilon$ iterations. In an iteration, the algorithm requires $O(n)$ computations at every vertex. This gives an $O(K_\epsilon n^2)$ running time for the BP algorithm to compute an $\epsilon$-approximate solution.

\section*{Acknowledgments}
This work resulted from discussions with Rajesh Sundaresan. I also thank him for suggestions that helped improve the presentation.
This work was supported by the Department of Science and Technology,
Government of India and by a TCS fellowship grant.


\bibliographystyle{imsart-number}
\bibliography{1}
\end{document}